\documentclass{amsart}

\usepackage{amssymb,amsmath,amsthm,graphics,amscd,amsfonts}

\theoremstyle{plain}
\newtheorem{theorem}{Theorem}[section]
\newtheorem{proposition}[theorem]{Proposition}

\newtheorem{example}[theorem]{Example}
\newtheorem{definition}[theorem]{Definition}

\newtheorem{remark}[theorem]{Remark}

\newcommand{\bfa}{{\bf a}}
\newcommand{\bfb}{{\bf b}}
\newcommand{\bfc}{{\bf c}}
\newcommand{\bfe}{{\bf e}}

\newcommand{\bfx}{{\bf x}}

\newcommand{\bfz}{{\bf z}}
\newcommand{\bfo}{{\bf o}}

\newcommand{\bfC}{{\mathbb C}}
\newcommand{\bfP}{{\mathbb P}}
\newcommand{\bfR}{{\mathbb R}}
\newcommand{\bfZ}{{\mathbb Z}}

\newcommand{\bfQ}{{\mathbb Q}}

\newcommand{\barj}{{\overline j}}

\newcommand{\barpartial}{{\overline \partial}}

\newcommand{\barz}{{\overline z}}

\newcommand{\tildeX}{\widetilde {X}}

\newcommand{\mapright}[1]{\smash{\mathop{   \hbox to 0.7cm{\rightarrowfill}}
 \limits^{#1}}}

\newcommand{\Ker}{{\rm Ker}}

\def\p{\partial}
\def\wt{\widetilde}

\def\p{\partial}

\def\a{\alpha}
\def\b{\beta}

\def\C{\mathbb C}

\def\ba{\begin{array}}
\def\ea{\end{array}}

\newtheorem{prop}[theorem]{Proposition}

\def\vol{\mathop{\mathrm{Vol}}\nolimits}
\def\grad{\mathop{\mathrm{grad}}\nolimits}

\begin{document}

\title
{Hilbert series and obstructions to asymptotic semistability}

\author{Akito Futaki}
\address{Department of Mathematics, Tokyo Institute of Technology, 2-12-1,
O-okayama, Meguro, Tokyo 152-8551, Japan}
\email{futaki@math.titech.ac.jp}
\author{Hajime Ono}
\address{Department of Mathematics,
Faculty of Science and Technology,
Tokyo University of Science,
2641 Yamazaki, Noda,
Chiba 278-8510, Japan}
\email{ono\_hajime@ma.noda.tus.ac.jp}

\author{Yuji Sano}
\address{Department of Mathematics, Kyushu University,
6-10-1, Hakozaki, Higashiku, Fukuoka-city, Fukuoka 812-8581 Japan}
\email{sano@math.kyushu-u.ac.jp}

\date{December 31, 2009}

\begin{abstract} 
Given a polarized manifold there are obstructions for asymptotic Chow semistability described as integral
invariants which can be regarded as characters of the Lie algebra of holomorphic vector fields. 
In 
this paper we show that, on toric Fano manifolds, the linear span of those Lie algebra characters coincides with the derivatives of
the Laurent series of the Hilbert series.

\end{abstract}

\thanks{The third author is supported by World Premier International
Research Center Initiative (WPI Initiative), MEXT Japan.}


\maketitle

\section{Introduction}
Let $M$ be a compact complex manifold of dimension $m$.
A polarization $L \to M$ is an ample line bundle over $M$, i.e. a holomorphic line bundle such that the first Chern class $c_1(L)$ is
represented by a positive $(1,1)$-form. 
Then $c_1(L)$ can be considered as a K\"ahler class. 
In \cite{donaldson01} Donaldson proved that if a polarized manifold $(M, L)$ admits a constant
scalar curvature K\"ahler metric (cscK metric for short) in $c_1(L)$ and if the automorphism group $\mathrm{Aut}(M)$ of $M$ is discrete then 
$(M,L)$ is asymptotically Chow stable. This result was extended by Mabuchi \cite{mabuchi-c} when $\mathrm{Aut}(M)$ is not discrete. Namely, Mabuchi proved that
if the obstruction introduced in \cite{mabuchi-a} vanishes and
$(M, L)$ admits a cscK metric in $c_1(L)$ then $(M, L)$ is asymptotically Chow polystable. 
The obstruction introduced in \cite{mabuchi-a} is an obstruction for $(M,L)$ to be asymptotically Chow semistable. We will explain
this obstruction and the definitions of relevant
stability conditions in section 2. We warn the reader that our terminology is slightly different from Mabuchi's.

Mabuchi's obstruction was reformulated by the first author in \cite{futaki04-1} to be the vanishing of a collection of integral invariants. One of these
integral invariants 
coincides with an obstruction to the existence of cscK metric (see \cite{futaki83.1}, \cite{futaki83.2}, \cite{calabi85}).
This last obstruction to the existence of cscK metrics is defined as a Lie algebra character on the complex Lie algebra ${\mathfrak h}(M)$ of all holomorphic vector fields on $M$, 
which we denote by $f : {\mathfrak h}(M) \to \bfC$. 
To explain the collection of integral invariants which obstruct the asymptotic semistability
let ${\mathfrak h}_0(M)$ the subalgebra consisting of all  holomorphic vector fields $X \in {\mathfrak h}(M)$ which have non-empty zero set.
Choose any $X \in {\mathfrak h}_0(M)$. 
For any
K\"ahler form $\omega$ representing $c_1(L)$ there exists a complex valued smooth
function $u_X$ determined up to a constant such that
$$ i(X) \omega = - \barpartial u_X. $$
When $u_X$ is a real function the real part of $X$ is a Hamiltonian vector field, and even if $u_X$ is not real we call $u_X$ the Hamiltonian function for $X$ by the abuse of terminology.
The existence of $u_X$ for $X \in {\mathfrak h}_0(M)$ is a classically known, see e.g.  \cite{kobayashitrnsf}; a comprehensive proof can be found in \cite{lebrunsimanca93}.
We assume the normalization of $u_X$ is so chosen that
\begin{equation}\label{normalization}
\int_M u_X\,\omega^m = 0.
\end{equation}
Thus, the Lie subalgebra ${\mathfrak h}_0(M)$ consists of all holomorphic vector fields $X$ in ${\mathfrak h}(M)$ such that $X$ is written in the form
\begin{equation}\label{grad}
 X = ig^{i\barj}\frac{\partial u_X}{\partial \bar{z}^j}\frac{\partial}{\partial z^i}
 \end{equation}
where the K\"ahler form $\omega$ is given by
$$ \omega = i g_{i\barj} dz^i \wedge d\barz^j. $$
To give another interpretation of $\mathfrak h_0(M)$, let $\mathrm{Aut}(M)$ be the group of all automorphisms of $M$. Let $\mathrm{Aut}(L)$ be the group of all bundle automorphisms
of $L$. Then $\mathrm{Aut}(L)$ contains $\bfC^{\ast}$ as a subgroup which acts as fiber multiplications. We put $\mathrm{Aut}(M,L) := \mathrm{Aut}(L)/\bfC^{\ast}$.
Then any element of $\mathrm{Aut}(M,L)$ induces an automorphism of $M$, and $\mathrm{Aut}(M,L)$ is considered as a Lie subgroup of $\mathrm{Aut}(M)$. 
The Lie subalgebra in $\mathfrak h(M)$ corresponding to  $\mathrm{Aut}(M,L)$ is exactly $\mathfrak h_0(M)$. This last fact follows from the general fact that
giving a moment map $M \to \mathfrak h_0(M)^\ast$ corresponds to giving a lifting of infinitesimal action of $\mathfrak h_0(M)$ on $M$ to that on $L$. 
Good references for this general fact are \cite{donkro}, section 6.5, and \cite{futakimabuchi02}, but the reference \cite{futakimabuchi02} is more precise in that the ambiguity
of Hamiltonian functions up to constant is more carefully treated. 
The functions $u_X$ define a moment
map, and thus define a lifting of the infinitesimal action of $\mathfrak h_0(M)$ on $M$ to an infinitesimal action on $L$. Therefore $\mathfrak h_0(M)$ corresponds to the Lie subgroup 
$\mathrm{Aut}(M,L)$ in $\mathrm{Aut}(M)$. However we encounter the trouble coming from the ambiguity of constants of Hamiltonian functions. Different constants give different
liftings. When we discuss stability we have to have a subgroup of special linear group. 
This point can be overcome by using S. Zhang's result (Theorem \ref{zhang} in section 2).

Let $\nabla$ be a type $(1,0)$ connection of the holomorphic tangent bundle, that is a linear connection whose connection form $\theta$ is expressed as a type $(1,0)$-form
with respect to local holomorphic frames. This last condition is of course equivalent to saying that the $(0,1)$-part of $\nabla$ is equal to $\bar\partial$.
Denote by $\Theta = \bar\partial\theta$ its curvature form. For a holomorphic vector field $X$ we also put
$$ L(X) := \nabla_X -  L_X  $$
where $L_X$ and $\nabla_X$ respectively denote the Lie derivative and covariant derivative by $X$. 
It is easy to see that $L(X)$ defines a smooth section of the endomorphism bundle of the holomorphic tangent bundle. 
Let $\phi$ be a $GL(m,\bfC)$-invariant polynomial of degree $p$ on ${\mathfrak gl}(m, \bfC)$.
We define ${\mathcal F}_{\phi} : \mathfrak h_0(M) \to \bfC$ by
\begin{eqnarray}\label{invariant}
{\mathcal F}_{\phi}(X) &=& (m-p+1) \int_M \phi(\Theta) \wedge u_X\,\omega^{m-p}
\\ & & + \int_M \phi(L(X) + \Theta) \wedge \omega^{m-p+1}.\nonumber
\end{eqnarray}
It can be shown that ${\mathcal F}_{\phi}(X)$ is independent of the choices of $\omega$ and $\nabla$, 
see \cite{futaki04-1} and \cite{FOSugaku} for the detailed account on ${\mathcal F}_{\phi}$. 
Let $\mathrm{Td}^p$ be the $p$-th Todd polynomial which is a $GL(m,\bfC)$-invariant polynomial of degree $p$ on ${\mathfrak gl}(m, \bfC)$.
The reformulation by the first author \cite{futaki04-1} of Mabuchi's obstruction \cite{mabuchi-a} to asymptotic Chow semistability of $(M,L)$ is the vanishing of
${\mathcal F}_{\mathrm{Td}^p}$ for all $p = 1,\,\cdots,\,m$. Moreover it can be shown that ${\mathcal F}_{\mathrm{Td}^1}$ coincides with the obstruction $f|_{\mathfrak h_0(M)}$ 
to the existence of cscK metric up to the multiplication of a non-zero constant. 

Choosing a K\"ahler form $\omega$ in $c_1(L)$, we have the Levi-Civita connection
$\theta = g^{-1} \partial g$ and its 
curvature form by $\Theta = \barpartial \theta$ where $g$ denotes the K\"ahler metric of $\omega$.
In the definition of ${\mathcal F}_{\phi}(X) $ above we could have used $\nabla$ to be the Levi-Civita connection
with the connection form $\theta$. 
In this case, because of the torsion-freeness of the Levi-Civita connection, $L(X)$ can be expressed by 
$$ L(X) = \nabla X = \nabla_j X^i\ dz^j \otimes \frac{\partial}{\partial z^i}$$
regarded as a smooth section of $\mathrm{End}(T'M)$.

The paper \cite{futaki04-1} was also motivated by the work of Mabuchi and Nakagawa \cite{mabuchi-nakagawa} in which they claimed that the
obstruction $f$ to the existence of cscK metric is an obstruction to semistability though their proof contained incomplete arguments, see the Erratum in \cite{mabuchi-nakagawa}.
But in light of their paper it is an interesting question whether or not 
the other integral invariants ${\mathcal F}_{\mathrm{Td}^2},\,\cdots,\,{\mathcal F}_{\mathrm{Td}^m}$ are linearly dependent on ${\mathcal F}_{\mathrm{Td}^1}$.

In this paper we related these integral invariants ${\mathcal F}_{\mathrm{Td}^p}$ to the Hilbert series for toric Fano manifolds, which is the index character considered by
Martelli, Sparks and Yau \cite{MSY2}.
The main result of this paper is Theorem \ref{Ind8} which claims that on toric Fano manifolds, the linear span of ${\mathcal F}_{\mathrm{Td}^2},\,\cdots,\,{\mathcal F}_{\mathrm{Td}^m}$ restricted
to the Lie algebra of the algebraic torus coincides with the linear span of the derivatives of
the Laurent series of the Hilbert series.
As an application we 
see that there are 3-dimensional toric Fano manifolds such that those integral invariants span 2-dimension. 
Thus ${\mathcal F}_{\mathrm{Td}^2},\,\cdots,\,{\mathcal F}_{\mathrm{Td}^m}$ are not in general linearly dependent on ${\mathcal F}_{\mathrm{Td}^1}$.

A question we can not answer in this paper is whether or not there is a polarized manifold $(M,L)$ on which a cscK metric exists in $c_1(L)$ so that ${\mathcal F}_{\mathrm{Td}^1} = 0$
but on which ${\mathcal F}_{\mathrm{Td}^p} \ne 0$ for some $p = 2,\, \cdots,\,m$. If the answer is no the assumption on the obstruction in Mabuchi's result \cite{mabuchi-c} can be omitted.
Our computations show that the last question is closely related to a question raised by Batyrev and Selivanova \cite{batyrev-selivanova}: Is a 
toric Fano manifold with vanishing $f$ for the anticanonical class necessarily symmetric\footnote{After posting the first version of this paper on the arXiv the paper \cite{NillP} appeared and a seven dimensional example of non-symmetric toric K\"ahler-Einstein Fano manifold
was presented. It is further shown in \cite{OSY} that, for this example, ${\mathcal F}_{\mathrm{Td}^1} = 0$ and ${\mathcal F}_{\mathrm{Td}^p} \ne 0$ for $p = 2,\, \cdots,\,7$.} ? 
If the answer is yes then any toric Fano K\"ahler-Einstein manifold has vanishing
${\mathcal F}_{\mathrm{Td}^p}$ for $p = 1,\,\cdots,\,m$.
Recall that a toric Fano manifold $M$ is said to be symmetric if the trivial character is the only fixed point of the action of the Weyl group on the space of all algebraic characters of the maximal 
torus in  $\mathrm{Aut}(M)$.
Note that if a toric Fano manifold $M$ is symmetric then the character $f$ for the anticanonical class vanishes. Recall also that Batyrev and Selivanova \cite{batyrev-selivanova} proved that 
a toric Fano manifold $M$ admits a K\"ahler-Einstein metric if $M$ is symmetric, and that Wang and Zhu \cite{wang-zhu} improved the result of Batyrev and Selivanova to the effect
that a toric Fano manifold $M$ admits a K\"ahler-Einstein metric if the invariant $f$ vanishes for the anticanonical class.

This paper is organized as follows. In section $2$, we review the
definition of asymptotic Chow semistability, and then we give the result,
Theorem \ref{SL2}, proved in \cite{futaki04-1}. In section $3$, we prove
the main theorem of this paper, Theorem \ref{Ind8}. In section $4$,
we give a combinatorial formula \eqref{MSY}
for computing the Hilbert series $C({\bf x},\mathcal{C}^*)$ of the toric
diagrams corresponding to the canonical bundle of toric Fano manifolds.
In section $5$, we observe that there are toric Fano $3$-folds such that
${\mathcal F}_{\mathrm{Td}^2}$ and ${\mathcal F}_{\mathrm{Td}^3}$ are
linearly independent of ${\mathcal F}_{\mathrm{Td}^1}$  using Theorem
\ref{Ind8} and computation of Hilbert series. In particular, we do computer
calculation in such a way as follows;

$\bullet$ Input the combinatorial data of the moment polytope to the 
formula \eqref{MSY}.

$\bullet$ Reduce the fractions to a common denominator.

$\bullet$ Substitute ${\bf e}^{-t{\bf b}}=(e^{-b_1t},e^{-b_2t},e^{-b_3t},
e^{-4t})$ for ${\bf x}=(x_1,x_2,x_3,x_4)$.

$\bullet$ Differentiate $C(e^{-t{\bf b}},\mathcal{C}^*)$ by $b_1,b_2$ and
$b_3$ at $(b_1,b_2,b_3)=(0,0,0)$.\\
In section $6$, we calculate ${\mathcal F}_{\mathrm{Td}^p}$ of toric
Fano $3$-folds by using localization formula.
We can also confirm the result obtained in section $5$ by this way.
In section $7$, we remark that a part of our construction makes sense
in general compact Sasaki manifolds.

\section{Obstructions to Asymptotic Chow semistability}

In this section we review the obstructions to asymptotic Chow semistability. Though the full account of these obstructions has already appeared in \cite{futaki04-1} we will reproduce
the arguments of \cite{futaki04-1} for the reader's convenience. The result in \cite{futaki04-1} we want to use in this paper is stated as follows.
\begin{theorem}[\cite{futaki04-1}]\label{SL2}
Suppose that a polarized manifold $(M,L)$ is asymptotically Chow semistable. 
Then ${\mathcal F}_{\mathrm{Td}^p}$ vanishes for all $p = 1,\,\cdots,\,m$. 
Suppose further that $H^q(M, \mathcal O(L))$ vanishes for all $q > 0$  and let $X$ be a holomorphic vector field in $\mathfrak h_0(M)$ which generates an $S^1$-action on $M$.
There is then a lifting of the infinitesimal action of $X$ to $L$ such that the induced $S^1$-action on $L^k$ induces an action of a subgroup in $\mathrm{SL}(H^0(M, \mathcal O(L^k)))$
for all $k \ge 1$ at once.
\end{theorem}

Note that in the above theorems the $S^1$-action on $L$ might be a finite covering of the $S^1$-action on $M$. 
Before reproducing the proof of Theorem \ref{SL2} we collect basic definitions and well-known facts used in the proof. First of all we recall the definitions of Chow (semi)stability 
and asymptotic Chow (semi)stability. 
\begin{definition}[Chow stability]\label{Chow1}
Let $\Lambda \to M$ be a very ample line bundle over a compact complex manifold $M$.
Let 
$\Phi_{|\Lambda|} : M \to \bfP(V)$ be the Kodaira embedding defined
by using the sections of $\Lambda$ where $V = H^0(M,\mathcal O(\Lambda))^{\ast}$. 
Let $d$ be the degree of $\Phi_{|\Lambda|}(M)$ in $\bfP(V)$. An element of the product
$\bfP(V^*) \times \cdots \times \bfP(V^*)$ of $m+1$ copies of $\bfP(V^*)$
defines  $m+1$ hyperplanes $H_1,\,\cdots\,H_{m+1}$ in $\bfP(V)$.
The set of all $m+1$ hyperplanes such that
$H_1 \cap \cdots \cap H_{m+1} \cap \Phi_{|\Lambda|}(M)$ is non-empty defines a divisor in
$\bfP(V^*) \times \cdots \times \bfP(V^*)$. Since the degree of $\Phi_{|\Lambda|}(M)$ is
$d$ this divisor is defined by some
${\hat M}_\Lambda \in (\mathrm{Sym}^{d}(V))^{\otimes (m+1)}$. Of course
${\hat M}_\Lambda$ is determined up to constant. The point
$[{\hat M}_\Lambda] \in \bfP( (\mathrm{Sym}^{d}(V))^{\otimes (m+1)})$ is called the Chow point.
$M$ is said to be Chow polystable with respect to $\Lambda$ if the orbit of ${\hat M}_\Lambda$
in $(\mathrm{Sym}^{d}(V))^{\otimes (m+1)}$ under the action of $\mathrm{SL}(V)$ is closed.
$M$ is said to be Chow stable with respect to $\Lambda$ if $M$ is polystable and the stabilizer at ${\hat M}_\Lambda$
of the action of  $\mathrm{SL}(V)$ is finite.
$M$ is said to be Chow semistable with respect to $\Lambda$ if the closure of the orbit of ${\hat M}_\Lambda$
in $(\mathrm{Sym}^{d}(V))^{\otimes (m+1)}$ under the action of $\mathrm{SL}(V)$ does not contain $\bfo \in (\mathrm{Sym}^{d}(V))^{\otimes (m+1)}$.
 Hilbert-Mumford criterion says that, to check Chow stability, it is sufficient to check the stability condition for all one parameter subgroups ($\cong \bfC^\ast$) in $SL(H^0(M,\mathcal O(\Lambda)))$.
\end{definition}
\begin{definition}[Asymptotic Chow stability]\label{Chow2}
Let $L \to M$ be an ample line bundle. 
For a large positive integer $k$, $L^k$ is very ample. We apply Definition \ref{Chow1} by taking $\Lambda$ to be $L^k$. 
Put $V_k := H^0(M,\mathcal O(L^k))^*$, and let
$\Phi_{|L^k|} : M \to \bfP(V_k)$ be the Kodaira embedding defined
by using the sections
of $L^k$. Let $d_k$ be the degree of $\Phi_{|L^k|}(M)$ in $\bfP(V_k)$. An element of the product
$\bfP(V_k^*) \times \cdots \times \bfP(V_k^*)$ of $m+1$ copies of $\bfP(V_k^*)$
defines  $m+1$ hyperplanes $H_1,\,\cdots\,H_{m+1}$ in $\bfP(V_k)$.
The set of all $m+1$ hyperplanes such that
$H_1 \cap \cdots \cap H_{m+1} \cap \Phi_{|L^k|}(M)$ is non-empty defines a divisor in
$\bfP(V_k^*) \times \cdots \times \bfP(V_k^*)$. Since the degree of $\Phi_{|L^k|}(M)$ is
$d_k$ this divisor is defined by some
${\hat M}_k \in (\mathrm{Sym}^{d_k}(V_k))^{\otimes (m+1)}$, 
which is determined up to constant. The point
$[{\hat M}_k] \in \bfP( (\mathrm{Sym}^{d_k}(V_k))^{\otimes (m+1)})$ is called the Chow point.
$M$ is said to be Chow polystable with respect to $L^k$ if the orbit of ${\hat M}_k$
in $(\mathrm{Sym}^{d_k}(V_k))^{\otimes (m+1)}$ under the action of $\mathrm{SL}(V_k)$ is closed.
$M$ is said to be Chow stable with respect to $L^k$ if $M$ is polystable and the stabilizer at ${\hat M}_k$
of the action of  $\mathrm{SL}(V_k)$ is finite.
$M$ is said to be Chow semistable with respect to $L^k$ if the closure of the orbit of ${\hat M}_k$
in $(\mathrm{Sym}^{d_k}(V_k))^{\otimes (m+1)}$ under the action of $\mathrm{SL}(V_k)$ does not contain $\bfo \in (\mathrm{Sym}^{d_k}(V_k))^{\otimes (m+1)}$.
$M$ is said to be asymptotically Chow polystable (resp. stable or semistable) with respect to $L$ if there exists a
$k_0 > 0$ such that $M$ is polystable  (resp. stable or semistable) for all $k \ge k_0$.
\end{definition}

Next recall a theorem of Zhang:
\begin{theorem}[ S.~Zhang, Theorem 3.4 in \cite{zhang96}]\label{zhang}\ \ Let $\Lambda \to M$ be a very ample line bundle over a compact complex manifold $M$,
and let 
$\Phi_{|\Lambda|} : M \to \bfP(V)$ be the Kodaira embedding defined
by using the sections of $\Lambda$ where $V = H^0(M,\mathcal O(\Lambda))^{\ast}$ as in Definition \ref{Chow1}. 
We identify $V$ with $\bfC^{N+1}$ 
endowed with the standard Hermitian metric. Then there is a norm, called Chow norm and denoted by $\Vert \cdot \Vert_C$, on
$(Sym^dV)^{\otimes(m+1)}$ such that, for any one parameter subgroup $\sigma_t$ of $\mathrm{SL}(V)$,
we have
\begin{equation}\label{zhang2}
\frac d{dt} \log \Vert \sigma_t(\hat{M}_\Lambda)\Vert_C = \int_M \dot\varphi_t \ \sigma^*_t\omega^m_\mathrm{FS}.
\end{equation}
where $\varphi_t = \log \Vert\sigma_t \bfz\Vert/\Vert \bfz \Vert$ for $\bfz \in V - 
\{\bfo \}$ with 
$[\bfz] \in M \subset \bfP(V)$, and
$\omega_{\mathrm{FS}}$ denotes the Fubini-Study K\"ahler form restricted to $M$.
\end{theorem}
\noindent
Zhang proved this theorem using Deligne pairing, but direct proofs are given also by Phong and Sturm \cite{phongsturm03} and the third author \cite{sano03}. 
Zhang's theorem lays a bridge between Chow stability and lifting an infinitesimal action of $\mathfrak h_0(M)$ to an ample line bundle $\Lambda$
in such a way that it induces a subgroup of $SL(H^0(M,\mathcal O(\Lambda)))$
for the following reasons.
Suppose that $\sigma_t$ in Theorem \ref{zhang} preserves $M$ and induces an action
on $M$ generated by
a holomorphic vector field $X$ on $M$. Then $\dot\varphi_t$ is a Hamiltonian function for $X$ with respect to the K\"ahler form $\sigma^*_t\omega_{FS}$.
For a polarized manifold $(M, \Lambda)$ with very ample line bundle $\Lambda$, let $V = H^0(M,\mathcal O(\Lambda))^{\ast}$ and $M \to \bfP(V)$ be the Kodaira
embedding. If 
$(M,\Lambda)$ is Chow semistable then the left hand side of (\ref{zhang2}) has to be zero, for otherwise the Chow norm tends to zero along the orbit of the one parameter subgroup $\sigma_t$
and the closure of the orbit is the straight line joining $\hat M_\Lambda$ and the origin $\bfo$, and of course contains $\bfo$. 
It follows therefore that if $(M,\Lambda)$ is Chow semistable then the Hamiltonian function $\dot\varphi_t = u_X$ 
giving a lifting to $\Lambda$ inducing a one parameter subgroup of $\mathrm{SL}(H^0(M,\mathcal O(\Lambda))$ satisfies 
 the normalization (\ref{normalization}) with respect to the
K\"ahler form $\sigma^*_t\omega_{FS}$. 
(Note that since $\Phi_{|\Lambda|}^\ast \mathcal O_{\bfP(V)} = \Lambda$ the restriction of the Fubini-Study metric to $M$ represents $c_1(\Lambda)$.)
Recall, as discussed in section 1, that 
the infinitesimal action of $X \in \mathfrak h_0(M)$ on $M$ always lifts to $\Lambda$ and this lifting is uniquely determined by a normalization of the Hamiltonian function
(refer to \cite{futakimabuchi02}),
which in this case is $\dot\varphi_t $. Therefore if $(M,\Lambda)$ is Chow semistable then the normalization (\ref{normalization}) gives a lifting of the infinitesimal action of $X  \in \mathfrak h_0(M)$ inducing a subgroup of 
$\mathrm{SL}(H^0(M,\mathcal O(\Lambda))$.
Summarizing the arguments given in this paragraph we get the following.
\begin{proposition}\label{zhang3} Let $(M, \Lambda)$ be a very ample line bundle and suppose that $(M, \Lambda)$ is Chow semistable. 
 Suppose also that we have a $\bfC^\ast$-action on $M$ generated by 
$X  \in \mathfrak h_0(M)$. Then the additive constant of Hamiltonian function $u_X$ for $X$ inducing a lifting of infinitesimal action of $X$ on $\Lambda$ such that it generates a subgroup
of $SL(H^0(M, \mathcal O(\Lambda)))$ is determined by the normalization (\ref{normalization}).
\end{proposition}

The detail of the following remarks the example can be found in \cite{futakimabuchi02}.

\begin{remark}\label{lift2}
Recall that two different liftings of $\bfC^\ast$-action on $M$ to $\Lambda$ differ by fiber multiplications by $\bfC^{\ast}$.
This difference of liftings is reflected by a difference of an additive constant of Hamiltonian functions. These constants belongs to $2\pi \bfZ$ corresponding to the weight of the
fiber multiplication. 
\end{remark}

\begin{remark}\label{lift4} Let $\Lambda \to M$ be an ample line bundle. 
Let $\omega$ be a K\"ahler form representing $c_1(\Lambda)$, and 
let $\tilde\theta$ be the connection form form of $\Lambda$ such that the curvature form of $\tilde\theta$ is $-2\pi i \omega$.
For an element $X$ in the  the lattice of the Lie algebra of the maximal torus in $\mathrm{Aut}(M,L)$,
a natural choice of a Hamiltonian function for $X$ with respect to $\omega$ is $(i/2\pi)\tilde\theta(X)$. Then 
 the normalization $\int_M (i/2\pi)\tilde\theta(X) \omega^m$ is a rational
number. This number is an invariant 
related to equivariant cohomology with respect to the lifted action. 
With this choice of Hamiltonian functions 
the moment map image is an integral polytope, i.e. a polytope with integer vertices.
In particular,  $u_X$ satisfying the normalization (\ref{normalization})
gives a rational polytope as its moment map image. 
\end{remark}

\begin{remark}\label{lift3} Let $\Lambda \to M$ be an ample line bundle. 
Consider an action of $S^1$ on $M$ generated by a holomorphic vector field $X\in \frak h_0(M)$ with period $1$. 
Suppose that a K\"ahler form $\omega \in c_1(\Lambda)$ is given.
Then the lift of the infinitesimal action of $X$ to $\Lambda$ is given locally by
$$ X \mapsto X^\sharp = -2\pi i u_X\, z \frac{\partial}{\partial z} + X^h$$
where $z$ is the fiber coordinate and $X^h$ is the horizontal lift with respect to the connection whose curvature is the given K\"ahler form $\omega$.
From Remark \ref{lift4} one can see that,  for $u_X$ satisfying the normalization (\ref{normalization}), $X^\sharp$ generates  an $S^1$-action
of integer period. See the next example.
\end{remark}

\begin{example} Let M be $\bfC\bfP^1$ and $ \Lambda$ be  $\mathcal O(1)$. Then the moment map image with respect to
the Hamiltonian using a cennection form is $[n, n+1] $ for some integer $n$. If we normalize the Hamiltonian function by (\ref{normalization}) then
the moment map image becomes $[-1/2, 1/2]$, and the period in this case is $2$.
\end{example}

Now we are in a position to give proofs of Theorem \ref{SL2}. To put it simply, applying Proposition \ref{zhang3} by
taking $\Lambda$ to be $L^k$ for all large $k$ shows that asymptotic stability assures that the normalization (\ref{normalization}) for $L$ gives
the normalization (\ref{normalization}) for $L^k$ all at once. We give some more detail about this in what follows.

\begin{proof}[Proof of Theorem \ref{SL2}]
We take an $X \in {\mathfrak h}_0(M)$ such that the real part $\mathrm{Re}X$ 
of $X$ generates an $S^1$-action $M$ with period $1$. 
Let $\omega$ be a K\"ahler form representing $c_1(L)$, and $u_X$ be the Hamiltonian function with the normalization (\ref{normalization}).
Then by Remark \ref{lift3}, $u_X$ defines a lifting of $X$ to an infinitesimal action of $X^\sharp$ which generates an $S^1$-action with
integer period. This induces a $\bfC^\ast$-action on $L$ and also on $L^k$ for all $k$ naturally. Fixing $k$, the lifted $\bfC^\ast$-action on $L^k$
defines a subgroup of $GL(H^0(M, \mathcal O(L^k)))$-action but not necessarily $SL(H^0(M, \mathcal O(L^k)))$. We therefore divide the
action by $\det^{1/N_k}$ where $N_k = \dim H^0(M, \mathcal O(L^k))$ and get a subgroup of $SL(H^0(M, \mathcal O(L^k)))$. This action
induces action on $V_k := H^0(M,\mathcal O(L^k))^*$, and thus on $L^k$. This last lifting the original $\bfC^\ast$-action on $M$ to $L^k$ must correspond
to a choice of Hamiltonian function of $X$ for the K\"ahler form $k\omega \in c_1(L^k)$ by the general principle. This Hamiltonian function must be of the form
\begin{equation}\label{Hamlift1}
\tilde u_{X,k} = ku_X + c_k
\end{equation}
where $c_k$ is a constant. Then
\begin{equation}\label{Hamlift2}
u_{X,k} = u_X + \frac{c_k}k
\end{equation}
defines a lifting of the infinitesimal action of $X$ on $M$ to an infinitesimal action $X^{\sharp}_k$ on $L$. 
By the construction, this $X^{\sharp}_k$ generates a $\bfC^\ast$-action on $L$ inducing a subgroup of $SL(H^0(M, \mathcal O(L^k)))$.

Now recall that we assume that the polarized manifold $(M,L)$ is asymptotically Chow semistable. Therefore there is a positive  integer $k_0$ such that
for all $k \ge k_0$, $(M,L^k)$ is Chow semistable. But by Proposition \ref{zhang3}, $\tilde u_{X,k}$ satisfies the normalization (\ref{normalization}). 
Since we also chose $u_X$ to satisfy this normalization we must have $c_k = 0$. This arguments apply for all $k \ge k_0$. 
Hence the choice of $u_X$ gives a lifting of $\bfC^\ast$-action to $L$ in such a way that the natural induced actions on $H^0(M, \mathcal O(L^k))$,
which we shall denote by $\rho_k$,  
defines a subgroup of $SL(H^0(M, \mathcal O(L^k)))$ for all $k \ge k_0$. 

Now we apply the equivariant Riemann-Roch theorem. 
Recall that we put $N_k = \dim H^0(M, \mathcal O(L^k))$. Then the weight $w_k$ of the action $\rho_k$
on $\wedge^{N_k}H^0(M, \mathcal O(L^k))$ is $0$ for $k \ge k_0$
since $\rho_k$ gives an $\mathrm{SL}$-action
on $H^0(M, \mathcal O(L^k))$. By the equivariant Riemann-Roch theorem
this weight is given by the coefficient of $t$ of the following (c.f. \cite{donaldson02}) : 
\begin{eqnarray*}
e^{k(\omega + t u_X)} \mathrm{Td}(tL(X) + \Theta)
= \sum_{p=0}^{\infty} \frac{k^p}{p!}(\omega + t u_X)^p
\sum_{q=0}^{\infty} \mathrm{Td}^{(q)}(tL(X) + \Theta).
\end{eqnarray*}
By writing the coefficient
of $t$ explicitly we have
\begin{equation}\label{equiRR}
0 = \sum_{p=0}^{m+1} \frac {k^p}{p!} \int_M (\omega^p \wedge 
\mathrm{Td}^{(m-p+1)}(L(X) + \Theta)  
+ p\,\omega^{p-1} \wedge u_X\,\mathrm{Td}^{(m-p+1)}(\Theta))
\end{equation}
for all $k \ge k_0$.
But from a result in \cite{futakimorita85} (see also Theorem 5.3.10 in \cite{futaki88})
\begin{equation}
\int_M \mathrm{Td}^{(m+1)}(L(X) + \Theta) = 0
\end{equation}
which implies that the term $p=0$ in (\ref{equiRR}) vanishes. The term $p = m+1$ also vanishes
because of our normalization (\ref{normalization}). 
Thus the vanishing of the terms for $p = 1,\ \cdots,\ m$ in (\ref{equiRR}) gives the
desired result since the terms for $p = 1,\ \cdots,\ m$ in (\ref{equiRR}) coincide with ${\mathcal F}_{\mathrm{Td}^p}$ for $p = 1,\,\cdots,\,m$. 

Conversely suppose that ${\mathcal F}_{\mathrm{Td}^p}$ vanishes for all $p = 1, \,\cdots,\,m$
and that $H^q(M, \mathcal O(L))$ vanishes for all $q > 0$, then the right hand side of (\ref{equiRR}) is zero
for any $k$. 
This implies that $\rho_k$ induces $\mathrm{SL}$-action not only for all $k \ge k_0$ but also for all $k \ge 1$. 
 This completes the proof of Theorem 
\ref{SL2}. 
\end{proof}

\begin{remark}\label{SL3}
Since $\mathrm{Td}^1 = \frac12 c_1$ and $c_1$ is the trace 
\begin{equation}\label{SL4}
{\mathcal F}_{\mathrm{Td}^1}(X) = \frac{m}2 \int_M S u_X \omega^m
\end{equation}
where $S$ denotes the scalar curvature of the K\"ahler form $\omega$. This is an obstruction to the existence of cscK metrics since if $S$ is constant
then ${\mathcal F}_{\mathrm{Td}^1}(X) = 0$ because of the normalization (\ref{normalization}). In fact the right hand side of (\ref{SL4}) is equal to $(mi/2)f(X)$.
To define $f$ 
let $F$ be a smooth function such that
$$ S - \int_M S \omega^m/\int_M \omega^m = \Delta F.$$
Then $f$ is define by 
$$ f(X) = \int_M XF\ \omega^m.$$
This $f(X)$ is independent of the choice of $\omega$ and obstructs the existence of a cscK metric in a given K\"ahler class  (\cite{futaki83.2}). 
Then using (\ref{grad}) we have 
\begin{eqnarray*}
\int_M S\, u_X\, \omega^m &=& \int_M \Delta F\, u_X\, \omega^m \\
&=& - \int_M (\grad u_X)F\ \omega^m \\
&=& i \int_M XF\ \omega^m = if(X).
\end{eqnarray*}
Therefore we get ${\mathcal F}_{\mathrm{Td}^1} = (mi/2)f(X)$.
\end{remark}
\begin{remark}\label{SL5}
If $M$ is a Fano and $L = K_M^{-1}$ 
it is more convenient to choose $F$ to be 
$$ \rho_{\omega} - \omega = \frac i{2\pi} \partial \barpartial F$$
where $\omega$ represents $c_1(M) = c_1(K_M^{-1})$,
the metric $g$ is given as 
$$ \omega = \frac i{2\pi} g_{i\barj} dz^i \wedge d\barz^j$$
and the Ricci form $\rho_{\omega}$ is given as 
$$ \rho_{\omega} = - \frac i{2\pi} \partial \barpartial \log \det (g_{i\barj}).$$
Then $f$ is defined as 
$$ f(X) = \int_M XF\, \omega^m.$$
Consider the second order elliptic differential operator
$$ \Delta_F = - g^{i\barj}\frac{\partial^2}{\partial z^i \partial \barz^j} - g^{i\barj}\frac{\partial F}{\partial z^i}\frac{\partial}{\partial \barz^j}.$$
If a complex valued smooth function $\wt v_X$ satisfies
$$ \Delta_F \wt v_X = \wt v_X$$
and put
$$ X := g^{i\barj}\frac {\partial \wt v_X}{\partial \barz^j} \frac{\partial}{\partial z^i}.$$
then we have
$$ i(X) \omega = i\barpartial \wt v_X$$
and
$$ f(X) = \frac 1{(m+1)}{\mathcal F}_{c_1^{m+1}}(X) = \int_M \Delta \wt v_X \rho_\omega^m = \int_M \mathrm{div}X\, \rho_\omega^m.$$
The proof of this result can be found in (5.2.1) in \cite{futaki88}. 
\end{remark}
\begin{remark}\label{SL6}
Mabuchi \cite{mabuchi-a} states the obstruction to asymptotic Chow semistability by
$$ \rho_{k} = \rho_{k_0} $$
for all $k \ge k_0$.
\end{remark}

\section{The Hilbert series}

Let $M$ be a toric Fano manifold of complex dimension $m$ and $K_M$ its canonical line bundle so that
the real torus $T^m$ acts on $M$ and this lifts to an action on $K_M$ by the pull-back of differential forms.
This $T^m$-action together with the $S^1$-action by multiplication on the fiber gives a $T^{m+1}$-action on $K_M$ so that
$K_M$ is also toric. 
It is a standard fact that $K_M^{-1}$ is very ample and $H^q(M, \mathcal O(K_M^{-1}))$ vanishes for all $q > 0$, see \cite{Oda} and \cite{Fulton}. 
We wish to consider the formal sum
$$ L(g) = \sum_{k=0}^\infty  Tr(g|_{H^0(M, \mathcal O(K_M^{-k}))})$$
where $Tr(g|_{H^0(M, \mathcal O(K_M^{-k}))})$ denotes the trace of the induced action of $g \in T^{m+1}$ on $H^0(M, \mathcal O(K_M^{-k}))$,
and regard $L(g)$ as a function of $g$. We call $L(g)$ the index character (c.f. \cite{MSY2}). We may analytically continue $L(\bfx)$ to $\bfx \in T^{m+1}_{\bfC}$,
the algebraic torus. 

Let $S$ be the total space of 
the associated $U(1)$-bundle of $K_M$. Then $S$ is a $(2m+1)$-dimensional Sasaki manifold. 
Recall that an odd dimensional Riemannian manifold $(S, g)$ 
is a Sasaki manifold if its Riemannian cone $(C(S), \bar g)$ with $C(S) = S \times \bfR_+$ and $\bar g = dr^2 + r^2 g$ is a K\"ahler manifold. Here $r$ denotes the
standard coordinate on $\bfR_+$. 
In the present case $C(S)$ is biholomorphic to $K_M - \{\mathrm{zero\ section}\}$, and $S$ is an $S^1$-bundle over the Fano manifold $M$. In such a case
we say that $S$ is a regular Sasaki manifold.

Since $M$ is toric so is $C(S)$. If the convex polytope (i.e. the moment map image) of $M$ is given by 
$$ P^{\ast} := \{ w \in \bfR^m\ |\ v_j \cdot w \ge -1\}$$
where $v_j \in \bfZ^m$ generates a $1$-dimensional face of the fan
then the convex polytope of $C(S)$ is given by 
$$ \mathcal C^{\ast} := \{ y \in \bfR^{m+1}\ |\ \lambda_j \cdot y \ge 0 \}$$
where $\lambda_j = (v_j, 1) \in \bfZ^{m+1}$. 
We denote by $P$ the dual polytope of $P^{\ast}$, that is, $P$ is a convex polytope with vertices $v_j \in \bfZ^m$. 
The integral points in $\mathcal C^{\ast}$ correspond to the sections of $K_M^{- k}$ for some $k \ge 1$.
In fact a point $(u,k) \in \bfZ^{m+1} \cap \mathcal C^{\ast}$ with $u \in \bfZ^m$ and $k \in \bfZ$ of height $k$ corresponds to an element $H^0(M, \mathcal O(K_M^{-k}))$. 
For $\bfa \in \mathcal C^{\ast} \cap \bfZ^{m+1} $ and $\bfx \in T^{m+1}_{\bfC}$ we put
$$ \bfx^{\bfa} := x_1^{a_1} \cdots x_{m+1}^{a_{m+1}}.$$
If $\bfa = (u,k)$ and $\sigma_\bfa$ denotes the holomorphic section of $H^0(M, \mathcal O(K_M^{-k}))$ corresponding to $\bfa$
then the action of $\bfx \in T^{m+1}_{\bfC}$ is given by $\sigma_\bfa \mapsto \bfx^{\bfa}\sigma_\bfa $. 
We write $C(\bfx, \mathcal C^{\ast})$ for the index character $L(\bfx)$ for the toric Fano manifold corresponding to the cone $\mathcal C^{\ast}$.
Thus we have obtained
$$ C(\bfx, \mathcal C^{\ast}) = \sum_{\bfa \in \mathcal C^{\ast} \cap \bfZ^{m+1}} \bfx^{\bfa}.$$
The right hand side is also called the Hilbert series.
It is known that the Hilbert series
$C({\bf x},\mathcal{C}^{\ast})$ of a rational cone $\mathcal{C}^{\ast}$
can be written as a rational function of ${\bf x}$, see the books
\cite{BR} or \cite{MS} for this subject.

For $\bfb \in \bfR^{m+1}$ we write
$$ \bfe^{-t\bfb} = (e^{-b_1t}, \cdots, e^{-b_m t}, e^{-b_{m+1} t})$$
and consider 
$$C(\bfe^{-t\bfb}, \mathcal C^{\ast}) = \sum_{\bfa \in \mathcal C^{\ast} \cap \bfZ^{m+1}} e^{-t\langle \bfa, \bfb \rangle} .$$
Then $C(\bfe^{-t\bfb}, \mathcal C^{\ast})$ is a meromorphic function of $t$.

We choose $\bfb \in \bfR^{m+1}$ from the subset
\begin{equation}\label{Ind0}
\mathcal C_R := \{ \bfb \in \bfR^{m+1}\ |\ \bfb = (b_1, \cdots, b_m, m+1),\ (b_1, \cdots, b_m) \in (m+1)P\}.
\end{equation}
The intrinsic meaning of the subset $\mathcal C_R$ can be explained in the context of toric Sasakian geometry as follows (c.f. \cite{MSY1},
\cite{MSY2}, and also \cite{FOW}, \cite{CFO}). 

We start with a general Sasaki manifold $S$ so that its Riemannian cone $C(S)$ described above is a K\"ahler manifold. When $C(S)$ is a toric K\"ahler
manifold we say that $S$ is a toric Sasaki manifold. 
$S$ is identified with $\{ r = 1\} \subset C(S)$. 
The Reeb vector field is a vector field $\xi = J(\partial/\partial r)$ on $S \cong \{ r = 1\} $ where
$J$ denotes the complex structure on $C(S)$. It extends to a vector field
on $C(S)$ given by $J(r\partial/\partial r)$, which we also call the Reeb vector field. 
The Reeb vector field is a Killing vector field both on $S$ and $C(S)$, and can be regarded as an element $\xi$ of the Lie algebra $\mathfrak t^{m+1}$ of the torus $T^{m+1}$.
When the cone $C(S)$ is $\bfQ$-Gorenstein as a toric variety it can be shown that there is an element $\gamma \in \mathfrak t^{m+1\ast}$ 
such that the Reeb vector field satisfies 
\begin{equation}\label{Ind2}
\langle \lambda_j, \gamma \rangle = -1
\end{equation}
and
\begin{equation}\label{Ind1}
\langle \gamma, \xi \rangle = - m - 1
\end{equation}
where $\lambda_j$'s in $\mathfrak t^{m+1}$ determine the moment cone $ \mathcal C^{\ast}$ of $C(S)$ by
$$ \mathcal C^{\ast} := \{ y \in \mathfrak t^{m+1\ast}
\ |\ \langle \lambda_j , y\rangle \ge 0 \}.$$
The smoothness of $C(S)$ implies that $\lambda_j$'s form a basis over $\bfZ$ along each $1$-dimensional face of $\mathcal C^{\ast}$.
Thus by (\ref{Ind2}), $\gamma$ is uniquely determined from the toric data of $C(S)$. 
If we vary the Sasakian structure by changing the Reeb vector field keeping the toric structure of $C(S)$, 
then, since $\gamma$ is not varied, the Reeb vector field $\xi$ has to obey the condition (\ref{Ind1}). Thus the deformation space of
Sasaki structures with fixed toric structure of the cone is given by

\begin{equation}\label{Ind3}
\mathcal C_R := \{\xi  \in \mathcal C \subset \mathfrak t^{m+1} \ |\ \langle \gamma, \xi \rangle = - m - 1 \}
\end{equation}
where $\mathcal C$ is the dual cone of $\mathcal C^{\ast}$. The tangent space $T_{\xi}C_R$ of the deformation space $C_R$ 
at $\xi$ is isomorphic to
$$ \{X  \in \ \mathfrak t^{m+1} \ |\ \langle \gamma, X \rangle = 0 \}. $$

For a regular Sasaki manifold we can take 
a basis of the lattice $\mathfrak t^{m+1}_\mathbb Z=\Ker \{\exp:
\mathfrak t^{m+1}\to T^{m+1}\}$ such that
$\gamma$ and the Reeb vector field $\xi_0$ are denoted as
$$\gamma = (0, \cdots, 0, -1)\ \ \text{and}\  \ \xi_0 = (0, \cdots, 0, m+1).$$
Then the deformation space (\ref{Ind3}) of Sasaki structures with fixed toric structure in this case coincides with the space (\ref{Ind0}).
The tangent space $T_{\xi_0}C_R$ of the deformation space $C_R$ 
at $\xi_0$ is isomorphic to
\begin{equation}\label{Ind4}
\{ \bfc \in \bfR^{m+1}\ |\ \bfc = (c_1, \cdots, c_m, 0)\}.
\end{equation}

The subspace given by (\ref{Ind4}) has another intrinsic meaning. Recall that the cone $C(S)$ for the regular Sasaki manifold $S$ is
$K_M$ minus the zero section for a Fano manifold $M$. In the toric case $M$ admits an action of the $m$-dimensional torus $T^m$, 
and together with the circle action of the fiber of $S \to M$, $S$ admits an action of the $(m+1)$-dimensional torus $T^{m+1}$. 
This $(m+1)$-dimensional torus action also gives the toric structure of $C(S)$. Let us consider the liftings of the action of $T^m$ on $M$ 
to $K_M$. A natural choice is given by the pull-back of differential forms since $K_M$ is the bundle of $(m,0)$-forms. 
Any other choice differs from the natural choice by the action
along the fibers of $S \to M$. The different choices of the liftings of the $T^m$-action on $M$ to $K_M$ can be described in two ways.

First of all, $C(S)$, which is isomorphic to $K_M$ minus the zero section, admits $T^{m+1}$-action. Any lifting of $T^m$-action on $M$ to $K_M$ is
given by a 
subgroup of $T^{m+1}$. At the Lie algebra level, this subgroup corresponds to a sub-lattice of rank $m$ in $\bfZ^{m+1} \subset \mathfrak t^{m+1}$. It spans a 
hyperplane in $\mathfrak t^{m+1}$. In this manner we can regard the hyperplane (\ref{Ind4}) as a lifting of $T^m$-action on $M$ to $K_M$ or $C(S)$.

Secondly, the difference of the liftings are described in terms of the normalization of Hamiltonian functions as follows.
Let $X$ be 
a holomorphic vector field on $M$ such that $X$ is the infinitesimal generator of the
action of an $S^1$ in $T^m$ and that $\exp(X) = 1$. Choose a lifting of the $S^1$-action on $M$ to $K_M^{-1}$ and 
let $\tildeX$ be its infinitesimal generator.
Then any other lift of the $S^1$-action 
is given by an infinitesimal generator of the form
$2\ell \pi iz\partial/ \partial z + \tildeX$ for some integer $\ell$ where $z$ denotes the coordinate of the fiber of $K_M^{-1} \to M$. 
Then if $\wt\theta$ is a connection form on the principal $\bfC^{\ast}$-bundle associated with $K_M^{-1}$ then 
$$ \frac1{2\pi i}\wt\theta(2\ell \pi iz\partial/ \partial z + \tildeX) = \frac1{2\pi i}\wt\theta(\tildeX) + \ell.$$
Let $\rho$ be  a $T^m$-invariant K\"ahler form representing $c_1(M) = - c_1(K_M)$. Then by the Calabi-Yau theorem \cite{yau77}
there is another 
$T^m$-invariant K\"ahler form $\omega$ representing $c_1(M) = - c_1(K_M)$ 
such that the Ricci form $\rho_{\omega}$ is equal to $\rho$.
Let 
$\wt\theta$ be the connection form on the principal $\bfC^{\ast}$-bundle associated with $K_M^{-1}$ of the Hermitian connection $\nabla$ induced from the Levi-Civita connection of $\omega$. 
Since $(i/2\pi)\wt\theta(\tildeX) + \mathrm{constant}$ is a Hamiltonian function of $X$ 
for the Ricci form $\rho_{\omega}$ considered as a symplectic form
and since the liftings of $T^m$-action to $K_M^{-1}$ and $K_M$
have the natural correspondence the above arguments explain that the difference of the liftings are described in terms of normalizations of Hamiltonian functions. 
To make this correspondence definitive we need to decide the Hamiltonian functions for the natural lifting by the pull-back of $(m,0)$-forms.

In the next section we shall consider the derivative of $C(\bfe^{-t\bfb}, \mathcal C^{\ast})$ at
$\bfb = \xi_0$ 
along a vector in the tangent space $T_{\xi_0}C_R$ described as (\ref{Ind4}). 
For that purpose we claim the following.

\begin{proposition}\label{Ind5}  Let $M$ be a Fano manifold and take $c_1(M)$ as a K\"ahler class. The following three liftings of $T^m$-action on $M$ to $K_M$ coincide.
Here the lifted action to $K_M$ naturally induces a lifted action to $K_M^{-1}$ and vice versa, and they are identified. 
\begin{enumerate}
\item[(a)] The action on $K_M$ defined by the pull-back of $(m,0)$-forms.
\item[(b)] The lifted action defined by the subspace (\ref{Ind4}).
\item[(c)] The lifted action to $K_M^{-1}$  defined by the normalization of the Hamiltonian function $v_X$ for $X \in \mathfrak t^m \otimes \bfC$ by
\begin{equation}\label{Ind6}
\int_M v_X\ \omega^m = \frac i{2\pi}f(X).
\end{equation}
Here $v_X$ is a Hamiltonian function of $X$ in the sense that $i(X)\omega = - \barpartial v_X$, and $f(X)$ is the one given in Remark \ref{SL5}.
\end{enumerate}
\end{proposition}

\begin{proof}  First we see that (a) and (c) coincide. 
As above let $\rho$ be  a $T^m$-invariant K\"ahler form representing $c_1(M) = - c_1(K_M)$. Then by the Calabi-Yau theorem \cite{yau77}
there is another 
$T^m$-invariant K\"ahler form $\omega$ representing $c_1(M) = - c_1(K_M)$ 
such that the Ricci form $\rho_{\omega}$ is equal to $\rho$.
Express the K\"ahler form 
$\omega = i g_{i\barj} dz^i \wedge d\barz^j$ as in section 1, and consider its Levi-Civita connection on the tangent bundle
and the induced connection on $K_M^{-1}$ and $K_M$. The pull back action of $T^m$ on $K_M$ is identified with
the usual push forward action on $K_M^{-1}$. Let $X$ be a holomorphic vector field whose real part belongs to the Lie algebra of $T^m$, and let $\tildeX$ be its lift
to $K_M^{-1}$ induced by the push forward action. 
It is easy to compute that for the connection form $\wt\theta$ on $K_M^{-1}$ we have 
$$\wt\theta(\tildeX) = \mathrm{div} X = \sum_{i=1}^m \nabla_i X^i.$$
From $\rho = \frac i{2\pi} \barpartial \theta$ we see that
$$ i(X) \rho = - \frac i{2\pi} \barpartial \wt\theta(\wt X)$$
and $v_X = \frac i{2\pi}\wt\theta(\wt X) =  \frac i{2\pi}\mathrm{div} X$ is the Hamitonian function. That (a) and (c) coincide follows from this.

The equivalence between (b) and (c) follows from the arguments given in the proof of Proposition 8.10 in \cite{FOW}. 
To explain these arguments we recall basic terminologies in Sasaki Geometry.
The Reeb vector field $\xi$ defines a flow which has a transverse K\"ahler structure. This means that the local orbit spaces are open K\"ahler manifolds
and that they are patched together isometrically on their overlaps. On these local orbit spaces we have K\"ahler forms which can be lifted to $S$ and
form a global two form $\omega^T$ called the transverse K\"ahler form. The Ricci forms on local orbit spaces also lifted to $S$ to form a global two form
$\rho^T$ called the transverse Ricci form. On local orbit spaces of the Reeb flow we have a $\bar\partial$ and $\partial$ operators, denoted by $\bar\partial_B$ and $\partial_B$. 
When the Sasaki manifold $S$ has a $\bfQ$-Gorenstein cone $C(S)$ there exists a smooth function $h$ such that
\begin{equation}\label{h1}
\rho^T - (2m+2) \omega^T = i\partial_B\bar\partial_B h. 
\end{equation}
This function $h$ is ``basic'' in the sense that locally it is obtained by lifting a function on the local orbit space. 
Note that the coefficient $(2m + 2)$ comes from the normalization of the Sasaki metric so that the length of Reeb vector field to be $1$.
With these terminologies in mind, 
it is proved in the proof of Proposition 8.10 in \cite{FOW} that, 
on the toric Sasaki manifold $S$ with $\bfQ$-Gorenstein cone $C(S)$, 
the tangent space to $C_R$ is equal to
\begin{equation}\label{Ind7}
\{X\in \mathfrak t^{m+1}\, |
\, \Delta_B^h \wt v_X=(2m+2)\wt v_X\},
\end{equation}
where 
$$\wt v_X = i((\bar\partial - \partial) \log r )(X) = -i X\log r,$$
$r$ being the coordinate on $\bfR_+$ in $C(S) = S \times \bfR_+$ and where
$$ \Delta^h = -  g_B^{i\barj}\frac{\partial^2}{\partial z^i \partial \barz^j} - g_B^{i\barj}\frac{\partial h}{\partial z^i}\frac{\partial}{\partial \barz^j},$$
$g_B$ being the transverse K\"ahler metric.
In the case of the regular Sasaki manifold $S$ over a Fano manifold $M$, the Reeb vector field
is induced by $\xi_0 = (0, \cdots, 0, m+1)$. The hyperplane given by (\ref{Ind4}) is equal to (\ref{Ind7}). But in this situation $(2m+2)\omega^T = \omega$ and $h = ((2m+2)/2\pi)F$.
Then the equation $\Delta_B^h \wt v=2(m+1)\wt v$ is equivalent to $\Delta_F \wt v = \wt v$ where
$$ \Delta_F = - g^{i\barj}\frac{\partial^2}{\partial z^i \partial \barz^j} - g^{i\barj}\frac{\partial F}{\partial z^i}\frac{\partial}{\partial \barz^j}.$$
Here, as in Remark \ref{SL5}, we take $\omega = (i/2\pi)g_{i\barj} dz^i \wedge d\barz^j$. 
Thus (\ref{Ind7}) implies $\Delta_F \wt v_X = \wt v_X$, and we have
$$
\int_M \wt v_X\ \omega^m = - \int_M  (\wt v_X)^i F_i\ \omega^m = - \int_M XF\ \omega^m = - f(X).
$$
But we see that $v_X = \frac{-i}{2\pi}  \wt v_X$ is the Hamiltonian function in the sense of (c) and it satsifies 
\begin{equation}
\int_M v_X \, \omega^m = \frac i{2\pi} f(X).
\end{equation}
This proves that (b) and (c) define the same lifting.
\end{proof}

Consider the derivatives of the coefficients of the Laurent series in $t$ of the meromorphic function $C(\bfe^{-t\bfb}, \mathcal C^{\ast})$ at $\bfb = \xi_0$ 
in the directions of vectors in the tangent space $T_{\xi_0}C_R$ described as (\ref{Ind4}).  Then those derivatives are characters of $\frak k \otimes \bfC$.

\begin{theorem}\label{Ind8}
 The linear span of those derivatives described as above coincides with the linear span
of ${\mathcal F}_{\mathrm{Td}^1},\,\cdots,\,{\mathcal F}_{\mathrm{Td}^m}$ restricted to $\frak k \otimes \bfC$. 
\end{theorem}

\begin{proof}
First of all, for a square matrix $A$ we have as a general formula in linear algebra
\begin{equation}\label{trace}
\frac{d}{ds}|_{s=0}\mathrm{tr}\, e^{sA} =  \frac{d}{ds}|_{s=0} \det e^{sA} = \mathrm{tr}\, A
\end{equation}
where $\mathrm{tr}$ denotes the trace. 
For a tangent vector $\bfc \in T_{\xi_0}C_R$ we consider the action of
$\bfe^{-t(\xi_0 + s \bfc)}$
on $H^0(M, K_M^{-k})$ and take the derivative with respect to $s$ at $s = 0$.
Since 
$$ \bfe^{-t\xi_0} = (1, \cdots , 1, e^{-(m+1)t}),$$
$$ \bfa = (a_1, \cdots , a_m, k)$$
and
$$ (\bfe^{-t\xi_0})^\bfa = e^{-t(m+1)k},$$
$\bfe^{-t\xi_0}$ acts on $H^0(M, K_M^{-k})$ as a scalar multiplication by $e^{-k(m+1)t}$, and 
we see from the general formula (\ref{trace}) that 
the derivative of $C(\bfe^{-t(\xi_0 + s \bfc)},
\mathcal C^{\ast})$ with respect to $s$ at $s = 0$ is the sum
$\sum_{k=1}^{\infty} e^{-k(m+1)}\tilde{w}_k$ where $\tilde{w}_k$ is the weight of
the lifted action described in Proposition \ref{Ind5}. By the equivariant index theorem each $\tilde{w}_k$ is given by 
\begin{equation}\label{equiRR2}
\tilde{w}_k = \sum_{p=0}^{m+1} \frac {k^p}{p!} \int_M (\omega^p \wedge 
\mathrm{Td}^{(m-p+1)}(L(X) + \Theta)  
+ p\,\omega^{p-1} \wedge v_X\,\mathrm{Td}^{(m-p+1)}(\Theta))
\end{equation}
where $v_X$ satisfies the 
the normalization (\ref{Ind6}). 
Recall that the Hamiltonian function
$u_X$ used in the definition of ${\mathcal F}_{\mathrm{Td}^k}(X)$ satisfies 
the normalization
$$ \int_M u_X\, \omega^m = 0. $$
Thus, $v_X = u_X - if(X)/2\pi \vol(M)$. Inserting this into the right hand side of (\ref{equiRR2}) one sees that our $\tilde{w}_k$ differs from 
$\sum_{p=1}^m (k^p/p!) {\mathcal F}_{\mathrm{Td}^k}(X)$ by a multiple of $f(X)$. But since 
$$f(X) = \frac {-2i}m {\mathcal F}_{\mathrm{Td}^1}(X)$$
we are done.
\end{proof}

\section{The formula of Martelli-Sparks-Yau}

The Hilbert series $C({\bf x},\mathcal{C}^*)$ of a toric diagram
$\mathcal{C}^*$,
which is the image of the moment map of a toric Calabi-Yau manifold,
is getting into the limelight in String theory, especially
AdS/CFT correspondence, for example
see \cite{BFHH,FHHZ,MSY2}.
Let 
$$\mathcal{C}^*=\{x\in \mathbb{R}^{m+1};\boldsymbol{\lambda}_i
\cdot x\ge 0,i=1,\cdots
,d\}$$
be an $(m+1)$-dimensional toric diagram of height $1$.
Here $\boldsymbol{\lambda}_i=
(v_i^1,\dots,v_i^m,1)\in \mathbb{Z}^{m+1}$
for each $i$. For a fixed ${\bf b}\in \mathcal{C}_R$, the Laurent expansion of
$C({\bf e}^{-t{\bf b}},\mathcal{C}^*)$ at $t=0$ is written as
\begin{equation}
C({\bf e}^{-t{\bf b}},\mathcal{C}^*)=\frac{C_{-m-1}({\bf b})}{t^{m+1}}+
\frac{C_{-m}({\bf b})}{t^m}+\frac{C_{-m+1}({\bf b})}{t^{m-1}}+\cdots.
\end{equation}
In \cite{MSY2}, Martelli, Sparks and Yau showed that
the coefficient of the leading order term $C_{-m-1}({\bf b})$ is a
constant multiple of the
volume of a Sasaki manifold whose Reeb vector field is generated by
${\bf b}$. Moreover they proved that if we think of ${\bf b}$ as
variables then
the first variation of $C_{-m-1}({\bf b})$ is equal to the Sasaki-Futaki
invariant\footnote{Strictly speaking, they proved this in the case when
${\bf b}$ is a rational vector. The general case was verified by the first two
authors and G. Wang  in \cite{FOW}.}.
Hence it is natural to ask what are the other coefficient
$C_i({\bf b})$ and its first variation for each $i$.
One of our motivations to write this article is
that we want to know the answer to this question.
As we saw in the previous section, 
when $\mathcal{C}^*$ corresponds to the canonical bundle of a toric
Fano manifold,
the first variations of $C_i$ at ${\bf b}=(0,\dots,0,m+1)$
are the linear combinations of the integral invariants $\mathcal{F}_{Td^l}.$

\begin{example}\label{ex1}
Let
$$\mathcal{C}^*=\{{\bf v}+x^1{\bf v}_1+\cdots+x^n{\bf v}_n;x^1,\dots,x^n\ge 0\}
\subset \mathbb{R}^n$$
be a rational simplicial cone, i.e. ${\bf v}_1,\dots,{\bf v}_n\in \mathbb{Z}^n$
and these are linearly independent in $\mathbb{R}^n$. Then by Theorem $3.5$
of \cite{BR}, the Hilbert series $C({\bf x},\mathcal{C}^*)$
of $\mathcal{C}^*$ is

\begin{equation}\label{BR}
C({\bf x},\mathcal{C}^*)
=\frac{\sigma_\Pi({\bf x})}{(1-{\bf x}^{{\bf v}_1})\cdots
(1-{\bf x}^{{\bf v}_n})},
\end{equation}
where $\Pi$ is the half-open parallelepiped
$$\Pi=\{{\bf v}+x^1{\bf v}_1+\cdots+x^n{\bf v}_n; 0\le x^1,\dots,x^n<1 \}$$
and 
$$\sigma_\Pi({\bf x})=\sum_{{\bf a}\in \Pi\cap \mathbb{Z}^n}{\bf x}^{\bfa}.$$

For example let 
$$\mathcal{C}^*:=\{a(1,1)+b(-1,1);a,b\ge 0\}\subset \mathbb{R}^2.$$
Then, 
$$C((x,y);\mathcal{C}^*)=\frac{1+y}{(1-xy)(1-x^{-1}y)}.$$
\end{example}

In the case when $\mathcal{C}^*$ is the toric diagram of height $1$
corresponding to the canonical bundle of a toric Fano manifold,
Martelli, Sparks and Yau \cite{MSY2} gave the formula to compute
$C({\bf x},\mathcal{C}^*)$ combinatorially.
For example, applying the formula to the Hilbert series in Example \ref{ex1},
we easily see that
\begin{align*}
\frac{1+y}{(1-xy)(1-x^{-1}y)}&=
\frac{1}{(1-xy)(1-x^{-1})}+\frac{1}{(1-x^{-1}y)(1-x)}\\
&=C((x,y),\mathcal{C}_1)+C((x,y),\mathcal{C}_2),
\end{align*}
where $\mathcal{C}_1=\{a(1,1)+b(-1,0);a,b\ge 0\},\mathcal{C}_2=\{
a(-1,1)+b(1,0);a,b\ge 0\}$.
To prove the formula, they formally applied the
Lefschetz fixed point formula to noncompact manifold
$K_M$, the total space of the canonical bundle of a toric Fano
manifold $M$. But we can verify the same formula using only combinatorial
argument as follows.
Let ${\bf v}_1,\dots,{\bf v}_d\in \mathbb{Z}^m$
be the vertices of a Fano polytope
$P_M\subset \mathbb{R}^m$. Equivalently, ${\bf v}_1,\dots,{\bf v}_d$
are the generators of $1$-dimensional cones of the fan of an
$m$-dimensional
toric Fano manifold $M$. If we set $\boldsymbol{\lambda}_j=({\bf v}_j,1)
\in \mathbb{Z}^{m+1}$, then we see that the cone
$$\mathcal{C}_M^*=\{{\bf x}\in \mathbb{R}^{m+1};
\boldsymbol{\lambda}_j\cdot {\bf x}\ge 0,j=1,\dots,d\}$$
is a toric diagram of height $1$ corresponding to the canonical bundle $K_M$
of $M$. We can also describe this cone $\mathcal{C}^*$ as
$$\mathcal{C}_M^*=\left\{\sum_{j=1}^{k}a^j\boldsymbol{\mu}_j;a^1,\dots,
a^k\ge 0,
\right\}$$
where $\boldsymbol{\mu}_j=({\bf w}_j,1)\in \mathbb{Z}^{m+1}$,
${\bf w}_1,\dots,{\bf w}_k$ is the vertices of the 
polar polytope
$$P^\circ_M=\{{\bf y}\in \mathbb{R}^m;{\bf v}_j\cdot{\bf y}
\ge -1,j=1,\dots,d\}.$$
Let ${\bf e}_{j,1},\dots,{\bf e}_{j,m}\in \mathbb{Z}^m$ denote the generators
of the edges emanating from a vertex ${\bf w}_j$.
Note here that ${\bf e}_{j,1},\dots,{\bf e}_{j,m}$ is a basis of $\mathbb{Z}^m$
for each $j$ since $P^\circ_M$ is a Delzant polytope.
Hence the Hilbert series of the cone
$$\mathcal{C}_{j,l}=\{l{\bf w}_j+x^1{\bf e}_{j,1}+\cdots+x^m{\bf e}_{j,m};
x^1,\dots,x^m\ge 0\}$$
in $\mathbb{R}^m$ is
$$C(\Tilde{\bf x},\mathcal{C}_{j,l})=
\frac{\Tilde{\bf x}^{l{\bf w}_j}}{(1-\Tilde{\bf x}^{{\bf e}_{j,1}})\cdots
(1-\Tilde{\bf x}^{{\bf e}_{j,m}})}
$$
by \eqref{BR}. Here $\Tilde{\bf x}=(x_1,\dots,x_m)$.
Then, by Brion's formula \cite{Brion}
\footnote{Brion's formula can be also obtained applying
Lefschetz fixed point formula to toric Fano manifolds, see \cite{I}.},
see also Theorem $9.7$ of \cite{BR} or
Theorem $12.13$ of \cite{MS}, we see that
$$
\sum_{\Tilde{\bf a}\in lP^\circ_M}
\Tilde{\bf x}^{\Tilde{\bf a}}=\sum_{j=1}^k
C(\Tilde{\bf x},\mathcal{C}_{j,l})
=\sum_{j=1}^k\Tilde{\bf x}^{l{\bf w}_j}
\prod_{b=1}^m
\frac{1}{(1-\Tilde{\bf x}^{{\bf e}_{j,b}})}.
$$
Therefore we have
\begin{equation}\label{MSY}
\begin{split}
C({\bf x},\mathcal{C}_M^*) &=
\sum_{l=0}^\infty \left\{\sum_{{\bf a}\in
\mathcal{C}^*_{M,\mathbb{Z}}\cap \{a_{m+1}=l\}}
{\bf x}^{\bf a}\right\}
=\sum_{l=0}^\infty \left\{\sum_{\Tilde{\bf a}\in lP^\circ_M}
\Tilde{\bf x}^{\Tilde{\bf a}}\right\}x_{m+1}^l\\
&=\sum_{l=0}^\infty \left\{\sum_{j=1}^k
\Tilde{\bf x}^{l{\bf w}_j}
\prod_{b=1}^m
\frac{1}{(1-\Tilde{\bf x}^{{\bf e}_{j,b}})}
\right\}x_{m+1}^l
=\sum_{l=0}^\infty 
\left\{\sum_{j=1}^k
{\bf x}^{l\boldsymbol{\mu}_j}
\prod_{b=1}^m
\frac{1}{(1-\Tilde{\bf x}^{{\bf e}_{j,b}})}
\right\}\\
&=\sum_{j=1}^k
\frac{1}{1-{\bf x}^{\boldsymbol{\mu}_j}}
\prod_{b=1}^m
\frac{1}{(1-\Tilde{\bf x}^{{\bf e}_{j,b}})}.
\end{split}
\end{equation}
Here $\mathcal{C}^*_{M,\mathbb{Z}}=\mathcal{C}_M^*\cap \mathbb{Z}^{m+1}$.
This is the formula given in \cite{MSY2}.
Note here that
$$
D_j({\bf x},\mathcal{C}_M^*):=
\frac{1}{1-{\bf x}^{\boldsymbol{\mu}_j}}
\prod_{b=1}^m
\frac{1}{(1-\Tilde{\bf x}^{{\bf e}_{j,b}})}
$$
diverges at ${\bf x}=(1,\dots,1,e^{-b^{m+1}t})={\bf e}
^{-t(0,\dots,0,b^{m+1})}$ for each $j$.
However if we reduce the fractions to a common denominator,
at least in the case of $m=1,2$, and $3$,
we see that there is a Laurent polynomial $K_{\mathcal{C}^*}({\bf x})$
such that
\begin{equation}\label{character}
C({\bf x},\mathcal{C}_M^*)=\frac{K_{\mathcal{C}_M^*}({\bf x})}
{(1-{\bf x}^{\boldsymbol{\mu}_1})\cdots(1-{\bf x}^{\boldsymbol{\mu}_k})}
\end{equation}
and $K_{\mathcal{C}_M^*}({\bf x})$
converges at
${\bf x}=(1,\dots,1,e^{-b^{m+1}t})$ when $b^{m+1},t\not=0$.
We do not know a general proof of this fact, but can check it using a computer
in each cases of $m = 1, 2$ and $3$.
For example, let $\mathcal{C}^*$ be the $2$-dimensional cone given in 
Example $4.1$.
Then $1/(1-xy)(1-x^{-1})$ and $1/(1-x^{-1}y)(1-x)$ diverge at
$x=1$. On the other hand
$$C((1,y),\mathcal{C}^*)=\frac{1+y}{(1-y)^2}.$$

To calculate the Hilbert series of the toric diagram $\mathcal{C}_M^*$
associated with a toric Fano manifold $M$,
we will use \eqref{character}.

\section{Examples}

In this section we give some combinatorial data and calculations
associated with toric
Fano threefolds. 
We used a computer algebra system Maxima\footnote{Maxima is available
from http://maxima.sourceforge.net/.}
for computing Hilbert series. Of course you can also utilize
other systems, for example, Maple, Mathematica and so on.
Since the expressions involved in the calculation are long we omit
them in this article.

The equivalence classes of toric Fano threefolds
(or $3$ dimensional Fano polytopes) are classified by Batyrev completely see
\cite{B0} or \cite{batyrev99}: There are $18$ equivalence classes
\begin{align*}
&\mathbb{C}P^3,\mathcal{B}_1,\mathcal{B}_2,\mathcal{B}_3,
\mathcal{B}_4=\mathbb{C}P^2\times \mathbb{C}P^1,
\mathcal{C}_1,\mathcal{C}_2,\mathcal{C}_3=\mathbb{C}P^1\times 
\mathbb{C}P^1\times \mathbb{C}P^1,\\
&\mathcal{C}_4,\mathcal{C}_5,\mathcal{D}_1,\mathcal{D}_2,
\mathcal{E}_1,\mathcal{E}_2,
\mathcal{E}_3,
\mathcal{E}_4,\mathcal{F}_1,\mathcal{F}_2,
\end{align*}
and
for each equivalence class, the vertices of Fano polytope are specified.
Here we use the same symbols as in \cite{batyrev99}
to represent toric Fano threefolds.
Hence we can compute the Hilbert series of the toric diagram
associated with the canonical bundles of Fano threefolds using
the formula in the previous section.

Let $M$ be a toric Fano threefold and $W$ the set of fixed point of 
the action of the Weyl group on 
the space of all algebraic  characters of the maximal torus in  
$\mathrm{Aut}(M)$.
Then we see that $\dim W=0,1,2$.

\subsection{The case when $\dim W=0$}

Let $M$ be a Fano threefold with $\dim W=0$, that is $M=\mathbb{C}P^3,
\mathbb{C}P^2\times \mathbb{C}P^1,\mathbb{C}P^1\times \mathbb{C}P^1\times
\mathbb{C}P^1,\mathcal{C}_5,\mathcal{F}_1$.
In such case the Futaki invariant
vanishes and by the result of Wang and Zhu, \cite{wang-zhu},
$M$ admits a K\"ahler-Einstein
metric. Moreover we see the following by calculation.

\begin{prop}
Let $M$ be a toric Fano threefold with $\dim W=0$.
Then
\begin{equation}
\frac{\partial H_M}{\partial a}(0,0,0;t)=
\frac{\partial H_M}{\partial b}(0,0,0;t)=
\frac{\partial H_M}{\partial c}(0,0,0;t)=
0.
\end{equation}
Here
$$H_M(a,b,c;t)=C((e^{-at},e^{-bt},e^{-ct},e^{-4t}),\mathcal{C}_M^*).$$
Therefore we see that the first variation of
$C_{i}(a,b,c)$ at $(a,b,c)=(0,0,0)$ vanishes for each
$i=-4,-3,-2,\dots$.
\end{prop}

\begin{example}
We give the combinatorial date when $M=\mathbb{C}P^2\times \mathbb{C}P^1$.

$\bullet $ The vertices of the Fano polytope $P_M$:
$$
\begin{pmatrix}
{\bf v}_1 & {\bf v}_2 & {\bf v}_3 & {\bf v}_4 & {\bf v}_5 
\end{pmatrix}
=\begin{pmatrix}
0 & 0 & 1 & -1 & 0\\
0 & 0 & 1 & 0 & -1\\
1 & -1 & 0 & 0 & 0
\end{pmatrix}
$$

$\bullet$ The vertices of the polar polytope $P^\circ _M$:
$$
\begin{pmatrix}
{\bf w}_1 & {\bf w}_2 & {\bf w}_3 & {\bf w}_4 & {\bf w}_5 & {\bf w}_6
\end{pmatrix}
=\begin{pmatrix}
1 & -2 & 1 & 1 & -2 & 1\\
-2 & 1 & 1 & -2 & 1 & 1\\
-1 & -1 & -1 & 1 & 1 & 1
\end{pmatrix}
$$

$\bullet$ The edges $\{{\bf e}_{j,1},{\bf e}_{j,2},{\bf e}_{j,3}\}$
emanating from ${\bf w}_j$:

\begin{align*}
&\begin{pmatrix}
{\bf e}_{1,1} & {\bf e}_{1,2} & {\bf e}_{1,3}
\end{pmatrix}
=
\begin{pmatrix}
-1 & 0 & 0\\
1 & 1 & 0\\
0 & 0 & 1
\end{pmatrix}
,\ \ 
\begin{pmatrix}
{\bf e}_{2,1} & {\bf e}_{2,2} & {\bf e}_{2,3}
\end{pmatrix}
=
\begin{pmatrix}
1 & 1 & 0\\
-1 & 0 & 0\\
0 & 0 & 1
\end{pmatrix}\\
&\begin{pmatrix}
{\bf e}_{3,1} & {\bf e}_{3,2} & {\bf e}_{3,3}
\end{pmatrix}
=
\begin{pmatrix}
-1 & 0 & 0\\
0 & -1 & 0\\
0 & 0 & 1
\end{pmatrix},\ \ 
\begin{pmatrix}
{\bf e}_{4,1} & {\bf e}_{4,2} & {\bf e}_{4,3}
\end{pmatrix}
=
\begin{pmatrix}
-1 & 0 & 0\\
1 & 1 & 0\\
0 & 0 & -1
\end{pmatrix}\\
&\begin{pmatrix}
{\bf e}_{5,1} & {\bf e}_{5,2} & {\bf e}_{5,3}
\end{pmatrix}
=
\begin{pmatrix}
1 & 1 & 0\\
-1 & 0 & 0\\
0 & 0 & -1
\end{pmatrix},\ \ 
\begin{pmatrix}
{\bf e}_{6,1} & {\bf e}_{6,2} & {\bf e}_{6,3}
\end{pmatrix}
=
\begin{pmatrix}
-1 & 0 & 0\\
0 & -1 & 0\\
0 & 0 & -1
\end{pmatrix}
\end{align*}
\end{example}

\subsection{The case when $\dim W=1$}

Let $M$ be a toric Fano threefold with $\dim W=1$. 
By the classification of toric Fano threefolds,
that is $M=\mathcal{B}_1,\mathcal{B}_2,\mathcal{B}_3,
\mathcal{C}_1,\mathcal{C}_4,\mathcal{E}_1,
\mathcal{E}_3,\mathcal{F}_2$.

\begin{prop}\label{2.3}
Let $M$ be a toric Fano threefold with $\dim W=1$.
Then
$$\left(\frac{\partial H_M}{\partial a}(0,0,0;t),
\frac{\partial H_M}{\partial b}(0,0,0;t),
\frac{\partial H_M}{\partial c}(0,0,0;t)
\right)=f(t){\bf p}.$$
Here ${\bf p}\in \mathbb{R}^3$ is a non-zero constant vector.
As a result the first variation of $C_i(a,b,c)$ at $(a,b,c)=(0,0,0)$
is a constant multiple
of that of $C_{-4}(a,b,c)$
for each $i=-3,-2,\dots$.
\end{prop}

\begin{example}
We give the combinatorial date when $M=\mathcal{B}_2$,
the blow-up of $\mathbb{C}P^3$ at a point.

$\bullet $ The vertices of the Fano polytope $P_M$:
$$
\begin{pmatrix}
{\bf v}_1 & {\bf v}_2 & {\bf v}_3 & {\bf v}_4 & {\bf v}_5 
\end{pmatrix}
=\begin{pmatrix}
1 & -1 & 0 & 0 & 1\\
0 & 0 & 1 & 0 & -1\\
0 & 0 & 0 & 1 & -1
\end{pmatrix}
$$

$\bullet$ The vertices of the polar polytope $P^\circ _M$:
$$
\begin{pmatrix}
{\bf w}_1 & {\bf w}_2 & {\bf w}_3 & {\bf w}_4 & {\bf w}_5 & {\bf w}_6
\end{pmatrix}
=\begin{pmatrix}
-1 & -1 & -1 & 1 & 1 & 1\\
-1 & -1 & 1 & -1 & -1 & 3\\
-1 & 1 & -1 & -1 & 3 & -1
\end{pmatrix}
$$

$\bullet$ The edges $\{{\bf e}_{j,1},{\bf e}_{j,2},{\bf e}_{j,3}\}$
emanating from ${\bf w}_j$:

\begin{align*}
&\begin{pmatrix}
{\bf e}_{1,1} & {\bf e}_{1,2} & {\bf e}_{1,3}
\end{pmatrix}
=
\begin{pmatrix}
1 & 0 & 0\\
0 & 1 & 0\\
0 & 0 & 1
\end{pmatrix}
,\ \ 
\begin{pmatrix}
{\bf e}_{2,1} & {\bf e}_{2,2} & {\bf e}_{2,3}
\end{pmatrix}
=
\begin{pmatrix}
1 & 0 & 0\\
0 & 1 & 0\\
1 & -1 & -1
\end{pmatrix}\\
&\begin{pmatrix}
{\bf e}_{3,1} & {\bf e}_{3,2} & {\bf e}_{3,3}
\end{pmatrix}
=
\begin{pmatrix}
1 & 0 & 0\\
1 & -1 & -1\\
0 & 1 & 0
\end{pmatrix},\ \ 
\begin{pmatrix}
{\bf e}_{4,1} & {\bf e}_{4,2} & {\bf e}_{4,3}
\end{pmatrix}
=
\begin{pmatrix}
-1 & 0 & 0\\
0 & 1 & 0\\
0 & 0 & 1
\end{pmatrix}\\
&\begin{pmatrix}
{\bf e}_{5,1} & {\bf e}_{5,2} & {\bf e}_{5,3}
\end{pmatrix}
=
\begin{pmatrix}
-1 & 0 & 0\\
0 & 1 & 0\\
-1 & -1 & -1
\end{pmatrix},\ \ 
\begin{pmatrix}
{\bf e}_{6,1} & {\bf e}_{6,2} & {\bf e}_{6,3}
\end{pmatrix}
=
\begin{pmatrix}
-1 & 0 & 0\\
-1 & -1 & -1\\
0 & 1 & 0
\end{pmatrix}
\end{align*}

$\bullet$ The ``gradient vector".
$$\left(\frac{\partial H_M}{\partial a}(0,0,0;t),
\frac{\partial H_M}{\partial b}(0,0,0;t),
\frac{\partial H_M}{\partial c}(0,0,0;t)
\right)=-\frac{te^{8t}(e^{4t}+3)(3e^{4t}+1)}{(e^{2t}-1)^5(e^{2t}+1)^5}
(3,1,1)$$
\end{example}

\subsection{The case when $\dim W=2$}

Let $M$ be a toric Fano threefold with $\dim W=2$, that is
$M=\mathcal{C}_2,\mathcal{D}_1,
\mathcal{D}_2,\mathcal{E}_2,\mathcal{E}_4$.
In this case, two different situations arise.

\begin{enumerate}
	\item When $M=\mathcal{C}_2$, then
	the ``gradient vector" is the same form as
	Proposition \ref{2.3}. Indeed we see that
	
	$\bullet $ The vertices of the Fano polytope $P_M$:
$$
\begin{pmatrix}
{\bf v}_1 & {\bf v}_2 & {\bf v}_3 & {\bf v}_4 & {\bf v}_5 & {\bf v}_6
\end{pmatrix}
=\begin{pmatrix}
0 & 0 & 1 & 1 & -1 & 0 \\
0 & 0 & 0 & 1 & 0 & -1\\
1 & -1 & 1 & 1 & 0 & 0
\end{pmatrix}
$$

$\bullet$ The vertices of the polar polytope $P^\circ _M$:
\begin{align*}
&\begin{pmatrix}
{\bf w}_1 & {\bf w}_2 & {\bf w}_3 & {\bf w}_4 & {\bf w}_5 & {\bf w}_6
& {\bf w}_7 & {\bf w}_8
\end{pmatrix}\\
&=\begin{pmatrix}
0 & 0 & 1 & 1 & -2 & -2 & 1 & 1\\
0 & 1 & -1 & 1 & 0 & 1 & -3 & 1\\
-1 & -1 & -1 & -1 & 1 & 1 & 1 & 1
\end{pmatrix}
\end{align*}

$\bullet$ The edges $\{{\bf e}_{j,1},{\bf e}_{j,2},{\bf e}_{j,3}\}$
emanating from ${\bf w}_j$:

\begin{align*}
&\begin{pmatrix}
{\bf e}_{1,1} & {\bf e}_{1,2} & {\bf e}_{1,3}
\end{pmatrix}
=
\begin{pmatrix}
1 & -1 & 0\\
-1 & 0 & 1\\
0 & 1 & 0
\end{pmatrix}
,\ \ 
\begin{pmatrix}
{\bf e}_{2,1} & {\bf e}_{2,2} & {\bf e}_{2,3}
\end{pmatrix}
=
\begin{pmatrix}
1 & 0 & -1\\
0 & -1 & 0\\
0 & 0 & 1
\end{pmatrix}\\
&\begin{pmatrix}
{\bf e}_{3,1} & {\bf e}_{3,2} & {\bf e}_{3,3}
\end{pmatrix}
=
\begin{pmatrix}
-1 & 0 & 0\\
1 & -1 & 1\\
0 & 1 & 0
\end{pmatrix},\ \ 
\begin{pmatrix}
{\bf e}_{4,1} & {\bf e}_{4,2} & {\bf e}_{4,3}
\end{pmatrix}
=
\begin{pmatrix}
-1 & 0 & 0\\
0 & -1 & 0\\
0 & 0 & 1
\end{pmatrix}\\
&\begin{pmatrix}
{\bf e}_{5,1} & {\bf e}_{5,2} & {\bf e}_{5,3}
\end{pmatrix}
=
\begin{pmatrix}
1 & 0 & 1\\
0 & 1 & -1\\
-1 & 0 & 0
\end{pmatrix},\ \ 
\begin{pmatrix}
{\bf e}_{6,1} & {\bf e}_{6,2} & {\bf e}_{6,3}
\end{pmatrix}
=
\begin{pmatrix}
1 & 0 & 1\\
0 & -1 & 0\\
0 & 0 & -1
\end{pmatrix}\\
&\begin{pmatrix}
{\bf e}_{7,1} & {\bf e}_{7,2} & {\bf e}_{7,3}
\end{pmatrix}
=
\begin{pmatrix}
-1 & 0 & 0\\
1 & 1 & 1\\
0 & 0 & -1
\end{pmatrix},\ \ 
\begin{pmatrix}
{\bf e}_{8,1} & {\bf e}_{8,2} & {\bf e}_{8,3}
\end{pmatrix}
=
\begin{pmatrix}
-1 & 0 & 0\\
0 & -1 & 0\\
0 & 0 & -1
\end{pmatrix}
\end{align*}

$\bullet$ The ``gradient vector".	
\begin{align*}
&\left(\frac{\partial H_M}{\partial a}(0,0,0;t),
\frac{\partial H_M}{\partial b}(0,0,0;t),
\frac{\partial H_M}{\partial c}(0,0,0;t)
\right)\\
&=-\frac{te^{8t}(e^{4t}+3)(3e^{4t}+1)}{(e^{2t}-1)^5(e^{2t}+1)^5}
(1,-2,3)
\end{align*}
	\item When $M=\mathcal{D}_1,\mathcal{D}_2,\mathcal{E}_2,\mathcal{E}_4$,
	the ``gradient vector" has components which are linearly independent
	as functions of $t$-variable.
	For example, let $M=\mathcal{D}_2$.
	Then we see that
	
	$\bullet $ The vertices of the Fano polytope $P_M$:
$$
\begin{pmatrix}
{\bf v}_1 & {\bf v}_2 & {\bf v}_3 & {\bf v}_4 & {\bf v}_5 & {\bf v}_6
\end{pmatrix}
=\begin{pmatrix}
0 & 1 & 0 & 1 & -1 & 0 \\
0 & 1 & 0 & 1 & 0 & -1\\
1 & 0 & -1 & 1 & 0 & 0
\end{pmatrix}
$$

$\bullet$ The vertices of the polar polytope $P^\circ _M$:
\begin{align*}
&\begin{pmatrix}
{\bf w}_1 & {\bf w}_2 & {\bf w}_3 & {\bf w}_4 & {\bf w}_5 & {\bf w}_6
& {\bf w}_7 & {\bf w}_8
\end{pmatrix}\\
&=\begin{pmatrix}
1 & 1 & -1 & 1 & 1 & -2 & 1 & -2\\
1 & -1 & 1 & 1 & -2 & 1 & -2 & 1\\
-1 & -1 & -1 & 1 & 1 & 1 & 0 & 0
\end{pmatrix}
\end{align*}

$\bullet$ The edges $\{{\bf e}_{j,1},{\bf e}_{j,2},{\bf e}_{j,3}\}$
emanating from ${\bf w}_j$:

\begin{align*}
&\begin{pmatrix}
{\bf e}_{1,1} & {\bf e}_{1,2} & {\bf e}_{1,3}
\end{pmatrix}
=
\begin{pmatrix}
-1 & 0 & 0\\
0 & -1 & 0\\
0 & 0 & 1
\end{pmatrix}
,\ \ 
\begin{pmatrix}
{\bf e}_{2,1} & {\bf e}_{2,2} & {\bf e}_{2,3}
\end{pmatrix}
=
\begin{pmatrix}
-1 & 0 & 0\\
1 & 1 & -1\\
0 & 0 & 1
\end{pmatrix}\\
&\begin{pmatrix}
{\bf e}_{3,1} & {\bf e}_{3,2} & {\bf e}_{3,3}
\end{pmatrix}
=
\begin{pmatrix}
1 & 1 & -1\\
0 & -1 & 0\\
0 & 0 & 1
\end{pmatrix},\ \ 
\begin{pmatrix}
{\bf e}_{4,1} & {\bf e}_{4,2} & {\bf e}_{4,3}
\end{pmatrix}
=
\begin{pmatrix}
-1 & 0 & 0\\
0 & -1 & 0\\
0 & 0 & -1
\end{pmatrix}\\
&\begin{pmatrix}
{\bf e}_{5,1} & {\bf e}_{5,2} & {\bf e}_{5,3}
\end{pmatrix}
=
\begin{pmatrix}
-1 & 0 & 0\\
1 & 1 & 0\\
0 & 0 & -1
\end{pmatrix},\ \ 
\begin{pmatrix}
{\bf e}_{6,1} & {\bf e}_{6,2} & {\bf e}_{6,3}
\end{pmatrix}
=
\begin{pmatrix}
1 & 1 & 0\\
0 & -1 & 0\\
0 & 0 & -1
\end{pmatrix}\\
&\begin{pmatrix}
{\bf e}_{7,1} & {\bf e}_{7,2} & {\bf e}_{7,3}
\end{pmatrix}
=
\begin{pmatrix}
-1 & 0 & 0\\
1 & 1 & 0\\
0 & -1 & 1
\end{pmatrix},\ \ 
\begin{pmatrix}
{\bf e}_{8,1} & {\bf e}_{8,2} & {\bf e}_{8,3}
\end{pmatrix}
=
\begin{pmatrix}
1 & 1 & 0\\
0 & -1 & 0\\
-1 & 0 & 1
\end{pmatrix}
\end{align*}

$\bullet$ The ``gradient vector".	
\begin{align*}
&\frac{\partial H_M}{\partial a}(0,0,0;t)
=\frac{\partial H_M}{\partial b}(0,0,0;t)=
-\frac{te^{8t}(2e^{8t}+7e^{4t}+2)}{(e^{4t}-1)^5(e^{4t}+1)^5}\\
&\frac{\partial H_M}{\partial c}(0,0,0;t)
=-\frac{te^{8t}(4e^{8t}+13e^{4t}+4)}{(e^{2t}-1)^5(e^{2t}+1)^5}
\end{align*}

$\bullet$ The first variation of $C_i(a,b,c)$ at $(a,b,c)=(0,0,0)$
for $i=-4,-2,-1$\footnote{Note here that $2C_{-4}(a,b,c)=C_{-3}(a,b,c)$ holds.
This is because $C_{-3}$ corresponds to the total scalar curvature of the
Sasaki manifold whose Reeb vector field is generated by $(a,b,c)$.
Therefore $C_{-3}$ is a constant multiple of $C_{-4}$.
See \cite{MSY2}}.
\begin{align*}
&\frac{\partial C_{-4}}{\partial a}(0,0,0)=
\frac{\partial C_{-4}}{\partial b}(0,0,0)=-\frac{11}{1024},\ \ 
\frac{\partial C_{-4}}{\partial c}(0,0,0)=-\frac{21}{1024}\\
&\frac{\partial C_{-2}}{\partial a}(0,0,0)=
\frac{\partial C_{-2}}{\partial b}(0,0,0)=-\frac{13}{768},\ \ 
\frac{\partial C_{-2}}{\partial c}(0,0,0)=-\frac{9}{256}\\
&\frac{\partial C_{-1}}{\partial a}(0,0,0)=
\frac{\partial C_{-1}}{\partial b}(0,0,0)=-\frac{1}{192},\ \ 
\frac{\partial C_{-1}}{\partial c}(0,0,0)=-\frac{1}{64}\\
\end{align*}

\end{enumerate}

\section{Direct computations of ${\mathcal F}_{\mathrm{Td}^p}$}
In this section, we shall check the results of the previous section by using the localization formula as in \cite{futakimorita85}.
As we saw in the previous section 
 $\{{\mathcal F}_{\mathrm{Td}^p}\}_{p=1,2,3}$ on toric Fano threefolds
 span at most dimension 2. 
We first show that this is generally
true for any toric Fano manifolds.
\begin{theorem}\label{dependence}
There is a universal linear dependence relation among $\{{\mathcal F}_{\mathrm{Td}^p}\}_{p=1,2,3}$ on any anticanonically
polarized toric Fano threefold $(M, K_M^{-1})$.
\end{theorem} 
\begin{proof}
We shall replace $\{{\mathcal F}_{\mathrm{Td}^p}\}_{p=1,2,3}$ by the invariants $\{{\mathcal G}_{\mathrm{Td}^p}\}_{p=1,2,3}$ with respect to the lifted action in  Proposition \ref{Ind5}.
Namely ${\mathcal G}_{\mathrm{Td}^p}$ is defined by the right hand side of
\eqref{invariant} with the normalization \eqref{Ind6}. 
Let $M$ be a Fano threefold and $\omega$ be a K\"ahler form in $c_1(M)$.
Let $\eta\in c_1(M)$ be another K\"ahler form whose Ricci form $\rho_{\eta}$ equals to $\omega$.
As in Proposition \ref{Ind5}, for $X\in {\mathfrak h}_0(M)$ let $u_X$ be the Hamiltonian function satisfying the normalization (\ref{Ind6}).
Recall that $u_X$ satisfies
\[
	\Delta_{F}u_X=\Delta_{\eta}u_X=u_X,
\]
where $\Delta_\eta$ is the complex Laplacian with respect to $\eta$.
Remark that the sign of $\Delta_\eta$ is consistent with $\Delta_F$ and $\Delta^h$.
Then we have
\begin{eqnarray}
	\nonumber
		12{\mathcal G}_{\mathrm{Td}^2}(X)
	&=&
		2\int_M (c^2_1+c_2)(\Theta_\eta)\wedge (\Delta_\eta u_{X}) \rho_{\eta}
		+\int_M (c^2_1+c_2)(L_\eta(X)+\Theta_\eta)\wedge \rho_{\eta}^{2}
	\\
	\nonumber
	&=&
		2\int_M (c^2_1+c_2)(\Theta_\eta)\wedge c_1(L_\eta(X)) c_1(\Theta_\eta)
	\\
	\nonumber
	&&	
		+\int_M (c^2_1+c_2)(L_\eta(X)+\Theta_\eta)\wedge c_1^2(\Theta_\eta)
	\\
	\nonumber
	&=&
		\int_M (c_1^4+c_1^2c_2)(L_\eta(X)+\Theta_\eta)
	\\
	\label{eq1_sano}
	&=&
		4{\mathcal F}_{\mathrm{Td}^1}(X)
		+\int_M (c_1^2c_2)(L_\eta(X)+\Theta_\eta).
\end{eqnarray}
Also we have
\begin{eqnarray}
	\nonumber
		24{\mathcal G}_{\mathrm{Td}^3}(X)
	&=&
		\int_M (c_1c_2)(\Theta_\eta) (\Delta_\eta u_X)
		+\int_M (c_1c_2)(L_\eta(X)+\Theta_\eta)\wedge \rho_{\eta}
	\\
	\nonumber
	&=&
		\int_M (c_1c_2)(\Theta_\eta) c_1(L_\eta(X)) 
		+\int_M (c_1c_2)(L_\eta(X)+\Theta_\eta)\wedge c_1(\Theta_\eta)
	\\
	\label{eq2_sano}
	&=&
		\int_M (c_1^2c_2)(L_\eta(X)+\Theta_\eta).
\end{eqnarray}
Since ${\mathcal G}_{\mathrm{Td}^p}(X)-{\mathcal F}_{\mathrm{Td}^p}(X)$ equals to a multiple of ${\mathcal F}_{\mathrm{Td}^1}(X)$ for all $p$ and $X$ as pointed out in the proof of Theorem \ref{Ind8}, the linear span of $\{{\mathcal F}_{\mathrm{Td}^p}\}_{p=1,2,3}$ equals to the one of $\{{\mathcal G}_{\mathrm{Td}^p}\}_{p=1,2,3}$.
From (\ref{eq1_sano}) and (\ref{eq2_sano}), $\{{\mathcal F}_{\mathrm{Td}^p}\}_{p=1,2,3}$ satisfies at least one linear dependence relation.
\end{proof}

This shows that 
the dimension of the span of $\{{\mathcal F}_{\mathrm{Td}^p}\}_{p=1,2,3}$ is not more than two for any Fano manifold with the polarization $L = K_M^{-1}$. 
To determine the dimension of the span of $\{{\mathcal F}_{\mathrm{Td}^p}\}_{p=1,2,3}$, it is therefore sufficient to investigate the linear independence between 
$\int_M (c_1^2c_2)(L_\eta(X)+\Theta_\eta)$ and $\int_M (c_1^4)(L_\eta(X)+\Theta_\eta)$.
Since both of them are kind of the integral invariants in \cite{futakimorita85}, we can apply the localization formula for them.
More precisely, if $X$ only has isolated zeroes, then
\begin{eqnarray}
	\label{local1}
		\int_M (c_1^4)(L(X)+\Theta)
	&=&
		\sum_{i}\frac{(\mbox{tr}(L(X_{p_i})))^4}{\det L(X_{p_i})},
	\\
	\label{local2}
		\int_M (c_1^2c_2)(L(X)+\Theta)
	&=&
		\sum_{i}\frac{(\mbox{tr}(L(X_{p_i})))^2\cdot
		 c_2(L(X_{p_i}))}
		{\det L(X_{p_i})},
\end{eqnarray}
where $\mbox{Zero}(X)=\{p_i\}_i \subset M$.
As for the localization formula, also see \cite{futaki88}.

Now we are in position to do calculations on examples. 
Firstly, let us compute $\mathcal{D}_2$.  We already saw in 5.3, (b) that 
for $\mathcal{D}_2$, $\{{\mathcal F}_{\mathrm{Td}^p}\}_{p=1,2,3}$ span a two dimensional vector space.
$\mathcal{D}_2$ is described as the blow up of $\mathbb{CP}^2\times \mathbb{CP}^1$ along $\mathbb{CP}^1 \times \{ \mbox{a point}\}$.
Let $[Z_0: Z_1: Z_2]$ be the homogeneous coordinate on $\mathbb{CP}^2$ and $[X_0: X_1]$ be the homogeneous coordinate on $\mathbb{CP}^1$.
Then let us consider the blow up of $\mathbb{CP}^2\times \mathbb{CP}^1$ along $\{([0: Z_1: Z_2], [1:0]) \mid [Z_1: Z_2] \in \mathbb{CP}^1\}$.
We denote the blow up by
\[
	\pi: \mathcal{D}_2 \to \mathbb{CP}^2\times \mathbb{CP}^1.
\]
Let $\sigma(t, \alpha, \beta, \gamma)$ be a flow on $\mathbb{CP}^2\times \mathbb{CP}^1=\{([Z_0: Z_1: Z_2], [X_0: X_1])\}$ defined by
\begin{equation}\label{flow}
\left(\begin{array}{ccccc}
e^{\alpha t} & 0 & 0 & 0 & 0 \\
0 & e^{\beta t} & 0 & 0 & 0 \\
0 & 0 & 1 & 0 & 0 \\
0 & 0 & 0 & e^{\gamma t} & 0 \\
0 & 0 & 0 & 0 & 1
\end{array}\right)
\end{equation}
where $t$ is a time parameter and $\alpha, \beta, \gamma \in \mathbb{R}$.
Remark that $\sigma$ transforms $\{([0: Z_1: Z_2], [1:0]) \mid [Z_1: Z_2] \in \mathbb{CP}^1\}$ into itself. So $\{\sigma\}$ can be lifted as a flow on $\mathcal{D}_2$.
We denote it by the same $\sigma$.
Remark that for generic $\alpha, \beta$ and $\gamma$ the set of all fixed points under the flow consists of the following eight isolated points.
\begin{eqnarray}
	\nonumber
	&&
	p_1:=\pi^{-1}(([1:0:0],[1:0])), 
	\\
	\nonumber
	&&
	p_2:=\pi^{-1}(([1:0:0],[0:1])), 
	\\
	\nonumber
	&&
	p_3:=\pi^{-1}(([0:1:0],[0:1])),
	\\
	\nonumber
	&&
	p_4:=\pi^{-1}(([0:0:1],[0:1]))
	\\
	\nonumber
	&&
	p_5:=(([0:1:0],[1:0]), 0),\,\,
	p_6:=(([0:1:0],[1:0]), \infty),
	\\
	\nonumber
	&&
	p_7:=(([0:0:1],[1:0]), 0),\,\,
	p_8:=(([0:0:1],[1:0]), \infty).
\end{eqnarray}
In above, $\mathcal{D}_2$ is regarded as the submanifold in $\mathcal{D}_2\times \mathbb{CP}^1=\mathcal{D}_2\times (\mathbb{C}\cup \{\infty\})$ with codimension one.
Let $X$ be the holomorphic vector field on $M$ associated with $\sigma$.
As for $p_1$ and $p_2$, we have
\[
	L(X)=\mbox{diag} (\beta-\alpha, -\alpha, \mp\gamma).
\]
As for $p_3$, we have
\[
	L(X)=
	\mbox{diag}( \alpha-\beta, -\beta, \gamma).
\]
As for $p_4$, we have
\[
	L(X)= \mbox{diag}(\alpha,\beta, \gamma).
\]
As for $p_5$, we have
\[
	L(X)=\mbox{diag}(\alpha-\beta,-\beta, -\alpha+\beta-\gamma).
\]
As for $p_6$, we have
\[
	L(X)=\mbox{diag}(-\beta,-\gamma,\alpha-\beta+\gamma).
\]
As for $p_7$, we have
\[
	L(X)=\mbox{diag}(\alpha, \beta, -\alpha-\gamma).
\]
As for $p_8$, we have
\[
	L(X)=\mbox{diag}(\beta, -\gamma, \alpha+\gamma).
\]
From (\ref{local1}) we find 
\[
	\int_M (c_1^4)(L(X)+\Theta)
	=-22\alpha+11\beta+21\gamma.
\]
Also from (\ref{local2}) we find
\[
	\int_M (c_1^2c_2)(L(X)+\Theta)
	=2(-2\alpha+\beta+3\gamma).
\]
Since the one is not proportional to the other, their span is two dimensional.

Next we shall compute $\mathcal{C}_2$. As we saw in 5.3, (a)  that $\mathcal{C}_2$ is 
 an example for which the invariants do span only  a one dimensional space although the set of fixed point of the action of 
 the Weyl group on the space of all algebraic  characters of the maximal torus in  $\mathrm{Aut}(M)$ is two dimensional.
There exists only single example among toric Fano threefolds, denoted by $\mathcal{C}_2$ in the list of Batyrev, 
which is $\mathbb{P}_{S_1} (\mathcal{O}\oplus \mathcal{O}(l))$ where $S_1$ is the blow up of $\mathbb{CP}^2$ at a point and
$l^2=1$ on $S_1$, i.e., the curve $l$ has the self-intersection $+1$.
Let $[Z_0:Z_1:Z_2]$ be the homogeneous coordinate on $\mathbb{CP}^2$ and $S_1$ be the blow up of $\mathbb{CP}^2$ at $[1:0:0]$.
Then, $S_1$ can be regarded as the submanifold of $\mathbb{CP}^2\times \mathbb{CP}^1$
\[
	\{
		([Z_0:Z_1:Z_2], [X_0:X_1]) \mid X_0Z_2-X_1Z_1=0
	\}.
\]
Let $[Y_0:Y_1]$ be the homogeneous coordinate on the fibers of $\mathcal{C}_2$.
Let $\sigma(t, \alpha, \beta, \gamma)$ be a flow on $\mathbb{CP}^2\times \mathbb{CP}^1=\{([Z_0: Z_1: Z_2], [Y_0: Y_1])\}$ defined as (\ref{flow}).
Since
\[
\left(\begin{array}{ccc}
e^{\alpha t} & 0 & 0 \\
0 & e^{\beta t} & 0  \\
0 & 0 & 1  
\end{array}\right)
\]
fixes $[1: 0: 0]$ in $\mathbb{CP}^2$, so $\{\sigma\}$ induces a flow on $\mathcal{C}_2$.
We denote it by the same $\sigma$ as the previous case.
For generic $\alpha, \beta$ and $\gamma$ the set of all fixed points under the action of $\sigma$ consists of the following eight isolated points.
\begin{eqnarray}
	\nonumber
	&&
	p_1:=([1:0:0],[1:0],[1:0]), 
		\,\,
	p_2:=([1:0:0],[1:0],[0:1]), 
	\\
	\nonumber
	&&
	p_3:=([1:0:0],[0:1],[1:0]),
		\,\,
	p_4:=([0:1:0],[0:1],[0:1]),
	\\
	\nonumber
	&&
	p_5:=([0:1:0],[1:0], [1:0]),\,\,
	p_6:=([0:1:0],[1:0], [0:1]),
	\\
	\nonumber
	&&
	p_7:=([0:0:1],[0:1], [1:0]),\,\,
	p_8:=([0:0:1],[0:1], [0:1]).
\end{eqnarray}
Let $X$ be the holomorphic vector field on $M$ associated with $\sigma$.
As for $p_1$ and $p_2$, we have
\[
	L(X)=\mbox{diag}(\beta-\alpha,-\beta,\mp\gamma).
\]
As for $p_3$ and $p_4$, we have
\[
	L(X)=\mbox{diag}(-\alpha,\beta,\mp\gamma).
\]
As for $p_5$ and $p_6$, we have
\[
	L(X)=\mbox{diag}(\alpha-\beta,-\beta,\pm(\alpha-\beta-\gamma)).
\]
As for $p_7$ and $p_8$, we have
\[
	L(X)=\mbox{diag}(\alpha,\beta,\pm(\alpha-\gamma)).
\]
Then we have
\[
		\int_M (c_1^4)(L(X)+\Theta)
	=
		4\int_M (c_1^2c_2)(L(X)+\Theta)
	=
	-16(-4\alpha+2\beta+3\gamma).
\]
The above equality implies our desired conclusion.

\section {Extension to general Sasaki manifolds}

In this section we remark that the invariants defined by (\ref{invariant}) extend to compact Sasaki manifolds.
As was explained in section 3 the Reeb vector field $\xi$ on a Sasaki manifold $S$ is defined as $J\frac{\partial}{\partial r}$. 
Let ${\mathcal F}_\xi$ be the {\it Reeb foliation} on $S$
generated by $\xi$. 
It is convenient to extend $\xi$ to a vector field $\widetilde{\xi} = J(r\frac{\partial}{\partial r})$ on $C(S)$.
It is well known that $\widetilde{\xi}-iJ\widetilde{\xi}=\widetilde{\xi}+ir\frac {\p}{\partial r}$ is a holomorphic vector
field on $C(S)$, and thus there is an action on $C(S)$ of the holomorphic flow
generated by $\widetilde{\xi}-iJ\widetilde{\xi}$. The collection of local orbit spaces of this action defines a
transversely holomorphic structure on the Reeb foliation ${\mathcal F}_\xi$ in the following sense.
There is an  open covering $\{U_\a\}_{\a\in A}$ of $S$ and submersions
$\pi_\a:U_\a\to V_\a\subset 
\C^m$ along the leaves of ${\mathcal F}_\xi$ such that when $U_\a\cap U_\b\neq \emptyset$ the transformation
$$\pi_\a\circ \pi_\b^{-1}:\pi_\b(U_\a\cap U_\b)\to \pi_\a(U_\a\cap U_\b)$$
is biholomorphic. We then have $\partial$ and $\bar\partial$ operators on each $V_\a$. They define well-defined 
operators, denoted by $\partial_B$ and $\bar\partial_B$, on the basic forms on $S$. Here a differential form $\psi$ on $S$ is
said to be basic if
$$ i(\xi) \psi = 0 \ \ \mathrm{and}\ \ \mathcal L_\xi \psi = 0.$$
We also put $d_B^c = \frac i2 (\bar{\partial}_B - \partial_B)$.

Let $G$ be a complex Lie group. 
We say that a principal $G$-bundle $P$ over $S$ is transversely holomorphic if the transition function
from $P|_{U_\b}$ to $P|_{U_\a}$ on the overlap 
$U_\a \cap U_\b$ is a holomorphic $G$-valued function on $\pi_\b(U_\a\cap U_\b)$ for any $\a$ and $\b$. 
A connection on $P$ is said to be a type $(1,0)$ connection if the connection form on $P|_{U_a}$ consists of type $(1,0)$ components
on $V_\a$ and $G$. For a  type $(1,0)$ connection on $P$ let $\Theta$ be its curvature 2-form. Then $\Theta$ does not have type $(0,2)$
components.

A typical such principal bundle is the frame bundle of the normal bundle $\nu({\mathcal F}_\xi)$ of the Reeb foliation ${\mathcal F}_\xi$ with $G = GL(m, \bfC)$. 
The Levi-Civita connections given by the transverse K\"ahler metric on local orbit spaces naturally define a global connection on  $\nu({\mathcal F}_\xi)$.
This is a typical example of type $(1,0)$ connection.

Let $\eta_0$ be the dual 1-form of $\xi$ on $S$. Then $\eta_0$ is a contact 1-form and $\frac 12 d\eta_0$ gives a transverse K\"ahler form.
Any other Sasaki structure compatible with the Reeb vector field $\xi$ is given by the deformation of $\eta_0$ into $\eta = \eta_0 + 2d_B^c\varphi$
for a basic function on $S$. 
This transformation induces the usual K\"ahler deformation in the transverse direction since it deforms $\frac12 d\eta_0$ into
$$ \frac12 d(\eta_0 + 2d_B^c\varphi) = \frac 12 d\eta_0 + d_Bd^c_B \varphi = \frac 12 d\eta_0 + i\partial_B\barpartial_B \varphi .$$
Let $\mathcal E$ be the set of all such contact forms $\eta = \eta_0 + 2d_B^c\varphi$. 

We pick an $\eta \in \mathcal E$ and fix it for the moment.
Let $\frak h_0$ be the Lie algebra of all holomorphic vector fields on $C(S)$ commuting with $\widetilde{\xi}-iJ\widetilde{\xi}$. Then a vector field in $\frak h_0$ defines naturally
a vector field on $S$. By the abuse of notation we also denoted by $\frak h_0$ the Lie algebra of all such vector fields on $S$. For such a vector field $X \in \frak h_0$ we put
$$ u_X = \eta(X) - \int_S \eta(X) \eta\wedge (d\eta)^m/\int_S  \eta\wedge (d\eta)^m.$$

Let $I^p(G)$ denote the set of all $G$-invariant polynomials of degree $p$ on $\mathfrak g$. For any $\phi \in I^p(G)$ we define $\mathcal F_\phi : \frak g \to \bfC$ by
\begin{eqnarray}\label{family}
{\mathcal F}_{\phi}(X) &=& (m-p+1) \int_S \phi(\Theta) \wedge u_X\,(d\eta)^{m-p}\wedge\eta
\nonumber
\\ & & + \int_S \phi(\theta(X) + \Theta) \wedge (d\eta)^{m-p+1}\wedge\eta.
\end{eqnarray}
Then one can prove the following theorem just as in \cite{futaki04-1} using Lemma 9.1 and Lemma 9.2 in \cite{FOW}.
\begin{theorem} ${\mathcal F}_{\phi}(X)$ is independent of the choices of $\eta \in \mathcal E$ and type $(1,0)$ connection $\theta$.
\end{theorem}
Let $L \to M$ be an ample line bundle. Then the total space $S$ of the associated $U(1)$-bundle is a Sasaki manifold.
If there is a torus action of the Sasaki structure we can deform the Sasaki structure by deforming the Reeb vector field , and 
we can consider $\mathcal F_\phi$ for irregular Sasaki manifolds. 

Let $M$ be a toric Fano manifold and  take $L$ to be $K_M^{-1}$. We can consider the integral invariants $\mathcal F_\phi$ for irregular Sasaki structures
obtained by deforming the Reeb vector field.
But it is not clear how 
the integral invariants $\mathcal F_{\mathrm{Td}^j}$ and the Hilbert series are related when the Reeb vector field is irregular. 

The following example is intriguing because it provides an example of a Sasaki manifold for which $\mathcal F_{\mathrm{Td}^1}$ vanishes but
the first variation of $C_{i}({\bf b})$
at the volume minimizing Reeb vector field
does not vanish for some $i$.
Let $M=\mathbb{C}P^2\#\overline{\mathbb{C}P^2}$. Then the total space of the associated $U(1)$-bundle has a Reeb field obtained by the 
volume minimization, and thus there is a Sasaki-Einstein metric. For this Reeb vector field $\mathcal F_{\mathrm{Td}^1}$ must vanish because
$\mathcal F_{\mathrm{Td}^1}$ is a multiple of the Sasaki-Futaki invariant. But the computation using a computer shows 
the first variation of $C_{-1}({\bf b})$
at this Reeb vector field
does not vanish.
We give the combinatorial data and some calculations
with respect to our $M$ below.

$\bullet $ The vertices of the Fano polytope $P_M$:
$$
\begin{pmatrix}
{\bf v}_1 & {\bf v}_2 & {\bf v}_3 & {\bf v}_4
\end{pmatrix}
=\begin{pmatrix}
1 & 0 & -1 & 0 \\
0 & 1 & -1 & -1 \\
\end{pmatrix}
$$

$\bullet$ The vertices of the polar polytope $P^\circ _M$:
$$
\begin{pmatrix}
{\bf w}_1 & {\bf w}_2 & {\bf w}_3 & {\bf w}_4 
\end{pmatrix}
=\begin{pmatrix}
2 & 0 & -1 & -1 \\
-1 & 1 & 1 & -1 
\end{pmatrix}
$$

$\bullet$ The edges $\{{\bf e}_{j,1},{\bf e}_{j,2}\}$
emanating from ${\bf w}_j$:

\begin{align*}
&\begin{pmatrix}
{\bf e}_{1,1} & {\bf e}_{1,2} 
\end{pmatrix}
=
\begin{pmatrix}
-1 & -1 \\
1 & 0
\end{pmatrix}
,\ \ 
\begin{pmatrix}
{\bf e}_{2,1} & {\bf e}_{2,2} 
\end{pmatrix}
=
\begin{pmatrix}
1 & -1 \\
-1 & 0 
\end{pmatrix}\\
&\begin{pmatrix}
{\bf e}_{3,1} & {\bf e}_{3,2} 
\end{pmatrix}
=
\begin{pmatrix}
1 & 0 \\
0 & -1 
\end{pmatrix},\ \ 
\begin{pmatrix}
{\bf e}_{4,1} & {\bf e}_{4,2} 
\end{pmatrix}
=
\begin{pmatrix}
1 & 0\\
0 & 1 
\end{pmatrix}\\
\end{align*}

$\bullet$ $C_i(a,b)$ ($i=-3,-2,-1$):

\begin{align*}
C_{-3}(a,b)&=\frac{2(2b-a+12)}{(b+3)(b-2a-3)(b-a+3)(b+a-3)},\ \ 
C_{-2}(a,b)=\frac32C_{-3}(a,b),\\
C_{-1}(a,b)&=-\frac{6b^2+2a^2b-6ab-18b-a^3+9a^2+9a-162}
{6(b+3)(b-2a-3)(b-a+3)(b+a-3)}
\end{align*}

$\bullet$ The volume minimizing $(a_0,b_0)$:
We call $(a_0,b_0)$ volume minimizing if the gradient of $C_{-3}(a,b)$
vanishes at $(a_0,b_0)$ and $(a_0,b_0)$ is in the interior of $3P^\circ _M$.
In this case, we see that
$(a_0,b_0)=(0,\sqrt{13}-4)$.
On the other hand,
$$
\frac{\partial C_{-1}}{\partial a}(a_0,b_0)=\frac{4(137\sqrt{13}-491)}
{(\sqrt{13}-7)^4(\sqrt{13}-1)^3},\ \ 
\frac{\partial C_{-1}}{\partial b}(a_0,b_0)=\frac{32(157\sqrt{13}-568)}
{(\sqrt{13}-7)^4(\sqrt{13}-1)^4}
$$
which is non-zero and hence the variation does not vanish.


\begin{thebibliography}{99}


\bibitem{B0}V.~V.~Batyrev, Toric Fano threefolds,
Izv. Skad. Nauk SSSR Ser. Mat. 45 (1981), no. 4, 704--717.

\bibitem{batyrev99} V.V. Batyrev, On the classification of toric Fano $4$-folds., J. Math. Sci. 94, (1999) 1021--1050.

\bibitem{batyrev-selivanova}
V.V. Batyrev and E.N. Selivanova : Einstein-K\"ahler metrics on symmetric toric Fano manifolds,  
J. Reine Angew. Math.  512 (1999) 225--236.



\bibitem{BR}M.~Beck and S.~Robins,
Computing the continuous discretely,
UTM, Springer, New York, 2007.

\bibitem{BFHH}S.~Benvenuti, B.~Feng, A.~Hanany and Y.-H.~He,
Counting BPS operators in gauge theories: quivers, syzygies and plethystics,
J. High Energy Phys. 2007, no. 11.

\bibitem{Brion}M.~Brion,
Points entiers dans les poly\`edres convexes,
Ann. Sci. \'Ecole Norm. Sup. (4), 21, 1988, 653--663.


\bibitem{calabi85}E.~Calabi : Extremal K\"ahler metrics II, Differential
geometry and complex analysis, (I. Chavel and H.M. Farkas eds.),  
95-114, Springer-Verlag, Berline-Heidelberg-New York,
(1985).

\bibitem{CFO} K.~Cho, A.~Futaki and H.~Ono : Uniqueness and examples of compact toric 
Sasaki-Einstein metrics, Comm. Math. Phys., 277 (2008), 439-458.  arXiv:math.DG/0701122.

\bibitem{donaldson01}S.K.~Donaldson : Scalar curvature and projective
embeddings, I, J. Differential Geometry, 59(2001), 479-522.

\bibitem{donaldson02}S.K.~Donaldson : Scalar curvature and stability of toric
varieties, J. Differential Geometry, 62(2002), 289-349.

\bibitem{donkro}S.K.~Donaldson and P.B.~Kronheimer : The geometry of four manifolds, Oxford Mathematical Monographs, Claren Press, Oxford, 1990.

\bibitem{FHHZ}D.~Forcella, A.~Hanany, Y.-H. He and A.~Zaffaroni,
The master space of $\mathcal{N}=1$ gauge theories,
arXiv:0801.1585.


\bibitem{Fulton}W.~Fulton : Introduction to Toric Varieties. Number 131, in Annals of Mathematics Studies,
Priceton University Press, Princeton, New Jersey, 1993. 

\bibitem{futaki83.1}A.~Futaki : 
An obstruction to the existence of Einstein K\"ahler metrics, Invent. 
Math. {\bf 73}, 437-443 (1983)

\bibitem{futaki83.2}A.~Futaki : On compact K\"ahler manifolds of constant scalar curvature, 
Proc. Japan Acad., Ser. A, {\bf 59}, 401-402 (1983)

\bibitem{futaki88}A.~Futaki : K\"ahler-Einstein metrics and integral 
invariants, Lecture Notes in Math., vol.1314, 
Springer-Verlag, Berline-Heidelberg-New York,
(1988)

\bibitem{futaki04-1}A.~Futaki : Asymptotic Chow semi-stability and integral
invariants, Intern. J. Math., {\bf 15}, 967-979, (2004).

\bibitem{futakimabuchi02}A.~Futaki and T.~Mabuchi : {\sl Moment maps and symmetric 
multilinear forms associated with
symplectic classes}, Asian J. Math.,
{\bf 6}, 349--372 (2002).


\bibitem{futakimorita85}A.~Futaki and S.~Morita : {\sl Invariant polynomials of the automorphism 
group of a compact complex manifold}, J. Differential Geometry,
{\bf 21}, 135--142 (1985).



\bibitem{FOSugaku}A.~Futaki and H.~Ono : Einstein metrics and GIT stability, 
to appear in Sugaku Expositions, arXiv:0811.0067. 
(Translation of Japanese version in Sugaku, 60(2008), 175--202).

\bibitem{FOW}A.~Futaki, H.~Ono and G.~Wang : Transverse K\"ahler geometry of
Sasaki manifolds and toric Sasaki-Einstein manifolds, to appear in J. Differential Geometry, math.DG/0607586.

\bibitem{I}A.~V.~S.~Infirri,
Lefschetz fixed-point theorem and lattice points in convex polytopes,
Adv. Math. 116 (1995), no. 1, 55--81.


\bibitem{kobayashitrnsf}S.~ Kobayashi, Transformation groups in differential geometry, 
Springer Verlag, Berlin-Heidelberg-New York, 1972.

\bibitem{lebrunsimanca93}C.~LeBrun and R.S.~Simanca : Extremal K\"ahler metrics and
complex deformation theory, Geom. Func. Analysis, 4 (1994) 298--336

\bibitem{mabuchi-a}T.~Mabuchi : An obstruction to asymptotic semistability and
approximate critical metrics, Osaka J. Math., 41(2004), 463-472. math.DG/0404210.

\bibitem{mabuchi-c}T.~Mabuchi : An energy-theoretic approach to the
Hitchin-Kobayashi correspondence for manifolds, I, Invent. Math. 159(2004), 225-243.

\bibitem{mabuchi-nakagawa}T.~Mabuchi and Y. Nakagawa : The Bando-Calabi-Futaki character as an obstruction to asymptotic semistability, 
Math. Ann., 324(2002), 187-193; Errutum, Math. Ann., 330(2004), 627-630.

\bibitem{MSY1}
D.~ Martelli, J.~ Sparks and  S.-T.~ Yau : The geometric dual of
a-maximisation for toric Sasaki-Einstein manifolds, 
Comm. Math. Phys., {\bf 268} (2006), 39-65. 
arXiv:hep-th/0503183.

\bibitem{MSY2} D.~Martelli, J.~ Sparks and  S.-T. ~Yau : Sasaki-Einstein manifolds and volume minimisation, 
Comm. Math. Phys., {\bf 280} (2008), 611-673. 
arXiv:hep-th/0603021.

\bibitem{MS}E.~Miller and B.~Sturmfels : 
Combinatorial commutative algebra,
GTM 227, Springer-Verlag, New York, 2005.

\bibitem{NillP}B.~Nill and A.~Paffenholz : 
Examples of non-symmetric K\"ahler-Einstein toric Fano manifolds, arXiv:0905.2054.

\bibitem{Oda}T.~Oda : Convex bodies and algebraic geometry, Springer-Verlag , Berline-Heidelberg-New York, 1988. 

\bibitem{OSY}H.~Ono, Y.~Sano and N.~Yotsutani : 
An example of asymptotically Chow unstable manifolds with constant scalar curvature, preprint, arXiv:0906.3836.


\bibitem{phongsturm03}D.H.~Phong and J.~Sturm : Stability, energy functionals,
and K\"ahler-Einstein metrics, Comm. Anal. Geom., vol. 11(2003), 563-597, 
arXiv: math.DG/0203254.

\bibitem{sano03}Y.~Sano : On stability-criterion of complete intersections, J. Geom. Anal., 14(2004), 533-544.

\bibitem{yau77}S.-T.~Yau : On Calabi's conjecture and some new results in
algebraic geometry, Proc. Nat. Acad. Sci. USA, {\bf 74}, 1798-1799 (1977)


\bibitem{wang-zhu}
X.J. Wang and X. Zhu : K\"ahler-Ricci solitons on toric manifolds with positive first Chern class, Adv. Math. 188, No.1 (2004) 87--103.


\bibitem{zhang96}S.~Zhang : Heights and reductions of semi-stable varieties,
Compositio Math. 104(1996), 77-105.



\end{thebibliography}
\end{document}